\documentclass[10pt,legno]{article}\textwidth 160mm\textheight 235mm
\UseRawInputEncoding
\usepackage{amsfonts}
\usepackage{amssymb}
\usepackage{graphicx}
\usepackage{amsmath}
\usepackage{amsthm}
\usepackage{tikz}
\usepackage{enumitem}
\usepackage{cancel}
\usepackage{todonotes}
\usepackage{cite}
\usepackage{hyperref}
\hypersetup{colorlinks=true,urlcolor=blue}
\usetikzlibrary{matrix}
\usetikzlibrary{arrows,shapes}
\usepackage{algorithm}
\expandafter\let\expandafter\oldproof\csname\string\proof\endcsname
\let\oldendproof\endproof
\renewenvironment{proof}[1][\proofname]{\oldproof[\ttfamily\scshape\bf #1.]}{\oldendproof}

\def\ve{\varepsilon}

\def\tilde{\widetilde}
\def\emp{\emptyset}

\def \dist {{\rm dist}}
\def\dom{{\rm dom}\,}
\def\d{{\rm d}}

\def\epi{{\rm epi\,}}

\def\N{{\cal N}}

\def\O{{\cal O}}

\def\d{{\rm d}}
\def\sub{\partial}

\def\B{\mathbb B}

\def\ox{\overline{x}}

\def\oy{\overline{y}}

\def\disp{\displaystyle}
\def\Lim{\mathop{{\rm Lim}}}

\def\Limsup{\mathop{{\rm Lim}\,{\rm sup}}}
\def\Liminf{\mathop{{\rm Lim}\,{\rm inf}}}

\def\tto{\;{\lower 1pt \hbox{$\rightarrow$}}\kern -10pt
\hbox{\raise 2pt \hbox{$\rightarrow$}}\;}
\def\Hat{\widehat}
\def\Tilde{\widetilde}
\def\Bar{\overline}
\def\ra{\rangle}
\def\la{\langle}
\def\ve{\varepsilon}
\def\B{I\!\!B}

\def\R{\mathbb{R}}
\def\N{\mathbb{N}}
\def\X{\mathbb{X}}
\def\Y{\mathbb{Y}}
\def\ox{\bar{x}}
\def\oy{\bar{y}}

\def\ov{\bar{v}}

\def\ou{\bar{u}}

\def\co{\mbox{\rm co}\,}
\def\cone{\mbox{\rm cone}\,}

\def\gph{\mbox{\rm gph}\,}
\def\epi{\mbox{\rm epi}\,}
\def\dim{\mbox{\rm dim}\,}
\def\dom{\mbox{\rm dom}\,}

\def\ker{\mbox{\rm ker}\,}

\def\dn{\downarrow}
\def\O{\Omega}
\def\ph{\varphi}
\def\emp{\emptyset}
\def\st{\stackrel}
\def\oR{\Bar{\R}}
\def\lm{\lambda}
\def\gg{\gamma}
\def\dd{\delta}
\def\al{\alpha}
\def\Th{\Theta}
\def\N{I\!\!N}
\def\th{\theta}
\def\vt{\vartheta}

\def\ss{\scriptsize }
\def\wto{\overset{w}{\rightarrow} }
\def\wsto{\overset{w^*}{\rightarrow} }

\def\wstoo{\overset{w^*}{\longrightarrow} }
\def\sce{\setcounter{equation}{0}}

\begin{document}
\newtheorem{Theorem}{Theorem}[section]
\newtheorem{Proposition}[Theorem]{Proposition}
\newtheorem{Remark}[Theorem]{Remark}
\newtheorem{Lemma}[Theorem]{Lemma}
\newtheorem{Corollary}[Theorem]{Corollary}
\newtheorem{Definition}[Theorem]{Definition}
\newtheorem{Example}[Theorem]{Example}
\newtheorem{Conclusion}[Theorem]{Conclusion}
\renewcommand{\theequation}{{\thesection}.\arabic{equation}}
\renewcommand{\thefootnote}{\fnsymbol{footnote}}
\begin{center}
{\bf\Large Variational Analysis in Normed Spaces\\with Applications to Constrained Optimization}\footnote{This research of was partly supported by the USA National Science Foundation under grants DMS-1512846, DMS-1808978 and by the USA Air Force Office of Scientific Research under grant \#15RT04.}\\[2ex] ASHKAN MOHAMMADI\footnote{Department of Mathematics, Wayne State University, Detroit, MI 48202, USA (ashkan.mohammadi@wayne.edu).}\;\mbox{ and }\;
BORIS S. MORDUKHOVICH\footnote{Department of Mathematics, Wayne State University, Detroit, MI 48202, USA (boris@math.wayne.edu). Research of this author was also partly supported by the Australian Research Council under Discovery Project DP-190100555.}
\end{center}

\vspace*{0.05in}
\small{\bf Abstract.} This paper is devoted to developing and applications of a generalized differential theory of variational analysis that allows us to work in incomplete normed spaces, without employing conventional variational techniques based on completeness and limiting procedures. The main attention is paid to generalized derivatives and subdifferentials of the Dini-Hadamard type with the usage of mild qualification conditions revolved around metric subregularity. In this way we develop calculus rules of generalized differentiation in normed spaces without imposing restrictive normal compactness assumptions and the like and then apply them to general problems of constrained optimization. Most of the obtained results are new even in finite dimensions. Finally, we derive refined necessary optimality conditions for nonconvex problems of semi-infinite and semidefinite programming.\\[0.5ex]
\vspace*{-0.2in}
\normalsize
\section{Introduction}\sce\vspace*{-0.05in}

A characteristic feature of modern variational analysis in both finite and infinite dimensions is a systematic usage of perturbation, approximation, and hence limiting procedures; see, e.g., the books \cite{bs,bz,c,i,m06,m18,penot,rw} and the references therein. In infinite dimensions these procedures require the {\em completeness} of the space in question that plays a crucial role in the development and applications of major {\em variational and extremal principles} of variational analysis. It has been well recognized that such principles (in particular, the most powerful Ekeland variational principle) are at the heart of variational techniques to derive calculus rules of generalized differentiation with their applications to optimization and control problems.

On the other hand, there exist---besides abstract optimization models---specific classes of infinite-dimensional optimization problems, which are naturally formulated in {\em incomplete} normed spaces. As important examples, we mention remarkable problems of the classical calculus of variations and optimal controls in the incomplete spaces ${\cal C}^1[a,b]$ and ${\cal C}^2[a,b]$ of continuously differentiable and twice continuously differentiable arcs, respectively.

The variational analysis developed in this paper in the normed space setting exploits generalized directional derivatives and subdifferentials of the {\em Dini-Hadamard type} for nonsmooth functions that are geometrically associated with the {\em tangent/contingent cone} to sets by Bouligand and Severi; see Section~\ref{prelim} for the definitions and discussions. These notions have been used in variational analysis and optimization theory in finite-dimensional and Banach space frameworks; see, e.g., \cite{bs,bn,i,ktz,m06,m18,penot,rw} and the references therein. However, in the majority of theoretical developments and applications they were considered as {\em preliminary/elementary} objects with the subsequent passage to the limit in order to create more robust constructions with adequate calculus rules based on variational and extremal principles in complete spaces. {\em Neither of such principles} is used in what follows.

Our novel approach here is quite different. We view the contingent/Dini-Hadamard constructions as the {\em major ones}, without their subsequent robust regularization, while examine them not for general functions and sets but under appropriate {\em qualification} and {\em regularity} conditions. For the case of compositions, the weakest qualification conditions used to derive the {\em equality-type chain rule} for {\em subderivatives} (or the Dini-Hadamard directional derivative) are the {\em Abadie} and {\em metric subregularity} ones. A special attention in the corresponding {\em subdifferential/normal cone} calculus for the Dini-Hadamard constructions is paid to a new class of {\em subamenable} functions and sets, which are also related to the metric subregularity condition that is much weaker that the more conventional metric regularity. It is shown, in particular, that the subamenability structure ensures the {\em robustness} (stability with respect to perturbations of the initial data) of the generalized differential constructions under consideration. While most of the obtained calculus results are new even in finite-dimensions, the developed calculus rules in infinite-dimensional normed spaces do {\em not} impose any {\em normal compactness} and the like assumptions, which are difficult to check. Recall that the latter assumptions play a crucial role in generalized differential calculus and applications of the major {\em limiting} subdifferentials, normal cones, and coderivatives in {\em Banach} spaces that are robust in much more general settings; see \cite{i,m06,penot} with the commentaries and bibliographies therein.\vspace*{0.03in}

In this paper we focus on developing generalized differential calculus in normed spaces with the subsequent applications to general problems of constrained optimization as well as to smooth nonlinear problems of semi-infinite and semidefinite programming in finite dimensions. Further applications to the aforementioned problems of the calculus of variations and optimal control in incomplete spaced will be given in a separate paper to follow. We also plan to extend these lines of research to second-order variational analysis and its applications in normed spaces.

The rest of this paper is organized as follows. Section~\ref{prelim} presents {\em basic definitions} and some preliminary results of variational analysis that are broadly employed in the formulations and proofs od the main results. In Section~\ref{sec03} we define the major {\em qualification conditions} used in the paper and establish relationships between them. Section~\ref{sec04} is devoted to obtaining calculus rules for {\em generalized directional derivatives} (subderivatives) under consideration with imposing the {\em weakest qualification conditions}. In Section~\ref{sec05} we turn to the corresponding {\em subdifferential calculus} in normed spaces with deriving the {\em equality-type} results under the refined notion of {\em Dini-Hadamard regularity}. Section~\ref{sec05a} mainly deals with {\em subamenable sets}, which are defined via the metric subregularity qualification condition, and develops rather surprising properties of the tangent and normal cones under consideration for this remarkable class of sets.

In Sections~\ref{sect06} and \ref{sect07} we present applications of the obtained calculus results to deriving {\em necessary optimality conditions} for some classes of constrained optimization problems. Section~\ref{sect06} concerns general class of optimization problems with the constraints described by $f(x)\in\Theta$ via a smooth mapping $f\colon\X\to\Y$ between {\em general normed spaces} and a closed set $\Theta\subset\Y$. We derive here two kinds of necessary optimality conditions for local minimizers: the {\em primal} and {\em dual} ones. The primal conditions are expressed in terms of the tangent cone to $\Th$ and the (directional) subderivative of a continuous cost function under the {\em Abadie constraint qualification}, which is a bit weaker than metric subregularity. The dual necessary optimality conditions are derived via the subdifferential of locally Lipschitzian and Dini-Hadamard regular cost functions and the normal cone to the convex set $\Th$ under the {\em metric subregularity constraint qualification}. In contrast to known results in this direction obtained in terms of some robust normal cones and subdifferentials, we add to the obtained multiplier rule a new {\em multiplier boundedness condition} with the bound expressed via the Lipschitz constant of the cost function multiplied by the constant taken from the imposed metric subregularity constraint qualification.

Section~\ref{sect07} is devoted to the application of the general necessary optimality conditions obtained in Section~\ref{sect06} to the two highly important classes that are well recognized in constrained optimization: {\em semi-infinite} and {\em semidefinite programs}in finite dimensions. The concluding Section~\ref{conc} summarizes the main achievements of the paper and discusses some direction of the future research.\vspace*{0.03in}

Throughout the entire paper we mostly use the standard notation of variational analysis and generalized differentiation; see, e.g., \cite{i,m06,rw}. Unless otherwise stated, all the spaces under consideration are arbitrary {\em normed spaces} endowed with the generic norm $\|\cdot\|$. Given such a space $\X$, the symbol $\la x,x^*\ra$ signifies the canonical pairing between $\X$ and it topological dual $\X^*$ with the indexes `$w$' and `$w^*$' signifying the weak and weak$^*$ topology on $\X$ and $\X^*$, respectively. For a nonempty set $\Omega\subset\X$, the notation $x\st{\O}{\to}\ox$ indicates that $x\to\ox$ with $x\in\O$, while $\co\O$ and $\cone\O$ stand for the convex and conic hulls of $\O$, respectively. The indicator function $\dd_\O$ of a set $\O$ is defined by $\dd_\O(x):=0$ for $x\in\O$ and $\dd_\O(x):=\infty$ otherwise, while the distance between $x\in\X$ and $\O$ is denoted by $\dist(x;\O)$. We write $x=o(t)$ with $x\in\X$ and $t\in\R_+$ indicating that $\frac{\|x\|}{t}\dn 0$ as $t\dn 0$, where $\R_+$ (resp.\ $\R_-$) means the collection of nonnegative (resp.\ nonpositive) real numbers. Recall also that $\N:=\{1,2,\ldots\}$.

We distinguish in notation between a single-valued $f\colon\X\to\Y$ and set-valued $F\colon\X\tto\Y$ mappings. Given a Fr\'echet differentiable mapping $f\colon\X\to\Y$ at $\ox$, its derivative (or Jacobian matrix in finite dimensions) is denoted by $\nabla f(\ox)$. In some results we also use the less restrictive notion of G\^ateaux differentiability and keep the same notation $\nabla f(\ox)$, while mentioning this in such cases. If $f\colon\X\to\Y$ is twice differentiable at $\ox$, the second derivative of $f$, which is a linear operator from $\X$ to the collections $L(\X,\Y)$ of linear bounded operators between $\X$ and $\Y$, can be identified with a continuous bilinear mapping $\nabla^2 f(\ox)\colon\X\times\X\to\Y$. If either $f$ is twice continuously differentiable around $\ox$, or both $\X$ and $\Y$ are finite-dimensional, then the bilinear mapping $\nabla^2f(\ox)(\cdot,\cdot)$ is symmetric, i.e., $\nabla^2f(\ox)(u,v)=\nabla^2f(\ox)(v,u)$ for all $u,v\in\X$. In the latter case we have by \cite[Theorem~13.2]{rw} the representation
\begin{equation*}
f(\ox+h)=f(\ox)+\la \nabla f(\ox),h\ra+\frac{1}{2}\nabla^2f(\ox)(h,h)+o(\|h\|^2).
\end{equation*}
If $\X=\R^n$ and $\Y=\R^m$ and if $f$ is twice differentiable at $\ox$, then
\begin{equation*}
\nabla^2f(\ox)(w,v)=\big(\la\nabla^2f_1(\ox)w,v\ra,\ldots,\la\nabla^2f_m(\ox)w,v\ra\big)\;\mbox{ for all }\;v,w\in\R^n,
\end{equation*}
where $f=(f_1,\ldots,f_m)$ and $\nabla^2f_i(\ox)$ stands for the Hessian matrix of $f_i$ at $\ox$.
\vspace*{-0.15in}

\section{Basic Definitions and Preliminaries}\label{prelim}\sce\vspace*{-0.05in}

We begin with recalling well-known tools of variational analysis and generalized differentiation utilized throughout the paper. The reader is referred to the books \cite{bz,c,i,m06,penot,rw} for more details and alternative terminologies in finite-dimensional and Banach spaces.

Given a parameterized family $(\O_{t})_{t>0}$ of nonempty subsets of $\X$, define the (Painlev\'e-Kuratowski) {\em outer limit} and {\em inner limit} of $\O_{t}$ as $t\dn 0$ by, respectively,
\begin{eqnarray*}
&&\Limsup\O_{t}:=\big\{x\in\X\;\big|\;\exists\,t_k\dn 0,\;\exists\,x_k\to x\;\mbox{ with }\;x_k\in\O_{t_k}\;\mbox{ as }\;k\in\N\big\},\\
&&\Liminf\O_{t}:=\big\{x\in\X\;\big|\;\forall\,t_k\dn 0,\;\exists\,x_k\to x,\;\exists\,k_0\in\N\;\mbox{ with }\;x_k\in\O_{t_k}\;\mbox{ as }\;k\ge k_0\big\}.
\end{eqnarray*}
The parameterized family of sets $\O_{t}$ {\em converges} to $\O$ as $t\dn 0$, i.e., $\O_t\to\O$, if
\begin{equation*}
\Limsup_{t\dn 0}\O_{t}=\Liminf_{t\dn 0}\O_{t}=\O:=\Lim_{t\dn 0}\O_{t}.
\end{equation*}
The (Bouligand-Severi) {\em tangent/contingent cone} $T_\O(\ox)$ to $\O$ at $\ox\in\O$ is defined by
\begin{equation}\label{tan}
T_\O(\ox):=\big\{u\in\X\;\big|\;\exists\,t_k\dn 0,\;\exists\,u_k\to u\;\mbox{ as }\;k\to\infty\;\mbox{ with }\;\ox+t_k u_k\in\O\big\}.
\end{equation}
We say that a tangent vector $u\in T_\O(\ox)$ is {\em derivable} if there exists $\xi\colon[0,\ve]\to\O$ with $\ve>0$, $\xi(0)=\ox$, and $\xi'_+(0)=u$, where $\xi'_+$ signifies the classical right derivative of $\xi$ at $0$. The set $\O$ is {\em geometrically derivable} at $\ox$ if every tangent vector $u$ to $\O$ at $\ox$ is derivable. The geometric derivability can be equivalently described by saying that the outer and inner limits of the parameterized family of sets $(\frac{\O-\ox}{t})_{t>0}$ agree as $t\dn 0$; in the other words, $\Lim_{t\dn 0}\frac{\O-\ox}{t}= T_{\O}(\ox)$.

Note that the contingent cone \eqref{tan} may be nonconvex in simple situations (e.g., for the graph $\O:=\gph|x|\subset\R^2$ at $\ox=(0,0)$), but the polar/dual to it defined by
\begin{equation}\label{dual}
N_{\O}^-(\ox):= T_{\O}(\ox)^*=\big\{v\in\X^*\;\big|\;\la v,u\ra\le 0\;\mbox{ for all }\;u\in T_{\O}(\ox)\big\}
\end{equation}
and known as the (Dini-Hadamard) {\em subnormal cone}, is always convex. This cone is generally larger than the (Fr\'echet) {\em regular normal cone}
\begin{equation}\label{fnc}
\Hat N_\O(\ox):=\disp\Big\{x^*\in\X^*\;\Big|\;\limsup_{x\st{\O}{\to}\ox}\frac{\la x^*,x-\ox\ra}{\|x-x_k\|}\le 0\Big\},
\end{equation}
which is also convex and can be equivalently described, if the space $\X$ is reflexive, in the dual scheme \eqref{dual} with the replacement of the contingent cone \eqref{tan} by its {\em weak contingent} counterpart $T^w_\O(\ox)$ defined in form \eqref{tan} while using the weak convergence $u_k\st{w}{\to}u$ instead of the strong one. Thus the cones \eqref{tan} and \eqref{fnc} agree in finite dimensions and the outer limit of them gives us in this case the (Mordukhovich) {\em limiting normal cone} that is defined in generality via the {\em sequential} outer limit
\begin{equation}\label{lnc}
N_\O^L(\ox):=\big\{x^*\in\X^*\;\big|\;\exists\,x_k\st{\O}{\to}\ox,\;\exists\,x^*_k\st{w^*}{\to}x^*\;\mbox{ with }\;x^*_k\in\Hat N_\O(x_k),\;k\in\N\big\}.
\end{equation}
Despite the nonconvexity of \eqref{lnc}, it enjoys---together with the associated constructions for nonsmooth functions and set-valued mappings---comprehensive calculus rules in the framework of Asplund spaces, i.e., such Banach spaces where each separable subspace has a separable dual. If $X$ is Asplund and $\O$ is closed, then the convex weak$^*$ closure of \eqref{lnc} agrees with the {\em convexified normal cone} $\Bar N_\O(\ox;\O)={\rm cl}^*{\rm co}\,N_\O^L(\ox)$ by Clarke, which is defined in Banach spaces by the duality scheme \eqref{dual} via his tangent cone to $\O$ at $\ox$.

If the set $\O$ is convex, then all the normal cones above reduce to the classical normal cone of convex analysis denoted by $N_\O(\ox)$, but for general nonconvex sets in normed spaces our major construction here is the subnormal one \eqref{dual}. We aim at developing strong and useful calculus rules for it and related subderivative and subdifferential notions for functions under appropriate qualification conditions without performing any limiting procedure as in \eqref{lnc}.\vspace*{0.03in}

Recall next the notion of a generalized directional derivative for functions, which is closely related to (in fact generated by) the contingent cone \eqref{tan}. Let $\ph\colon\X\to\oR:=(-\infty,\infty]$ be an extended-real-valued function on a normed space with the associated domain and epigraph sets
\begin{equation*}
\dom\ph:=\big\{x\in\X\;\big|\;\ph(x)<\infty\big\}\;\mbox{ and }\;\epi\ph:=\big\{(x,\al)\in\X\times\R\;\big|\;\al\ge\ph(x)\big\}.
\end{equation*}
Then $\ph$ is said to be {\em proper} if $\dom\ph\ne\emp$. The (Dini-Hadamard) {\em subderivative} of $\ph$ at $\ox\in\dom\ph$ is the function ${\mathrm d}\ph(\ox)\colon\X\to[-\infty,\infty]$ with the finite and infinite values $\pm\infty$ defined by
\begin{equation}\label{subder}
{\mathrm d}\ph(\ox)(\ou):=\liminf_{\substack{t\dn 0\\u\to\ou}}{\frac{\ph(\ox+tu)-\ph(\ox)}{t}},\quad\ou\in\X.
\end{equation}
The subderivative construction \eqref{subder} is geometrically determined by the contingent cone \eqref{tan} via the relationship $T_{\ss\epi\ph}(\ox,\ph(\ox))=\epi\d\ph(\ox)$. This clearly implies that the subderivative function is lower semicontinuous and positive homogeneous of degree 1.

Further, the function $\ph$ is called {\em epi-differentiable} at $\ox$ if the family of sets $\epi\Delta_t\ph(\ox)(\cdot)$ with
\begin{equation}\label{Delta}
\Delta_t\ph(\ox)(u):=\frac{\ph(\ox+tu)-\ph(\ox)}{t},\quad u\in\X,
\end{equation}
converges to $\epi\d\ph(\ox)$ in the above sense of set convergence. If in addition ${\mathrm d}\ph(\ox)$ is a proper function, then $\ph$ is {\em properly epi-differentiable} at $\ox$. It is easy to see that the epi-differentiability of $\ph$ at $\ox$ agrees with the geometric derivability of $\epi\ph$ at $(\ox,\ph(\ox))$. Observe also that if $\ph$ is epi-differentiable at $\ox$, then its domain $\dom\ph$ is geometrically derivable at $\ox$.\vspace*{0.03in}

Next we consider the subdifferential notions for extended-real-valued functions that are generated by the corresponding notions of normals to sets discussed above. Given $\ph\colon\X\to\oR$ on a normed space $\X$, the {\em Dini-Hadamard subdifferential} of $\ph$ at $\ox\in\dom\ph$ is
\begin{equation}\label{sub}
\sub^{-}\ph(\ox):=\big\{v\in\X^*\;\big|\;(v,-1)\in N_{\ss\epi\ph}^{-}\big(\ox,\ph(\ox)\big)\big\}
\end{equation}
via the subnormal cone \eqref{dual} to the epigraph of $\ph$. Note that \eqref{sub} is a standard geometric scheme to define subdifferentials via normal cones, and all the normal cones mentioned above generate the corresponding subdifferentials (collections of subgradients). On the other hand, we can get normal cones to sets from subdifferentials of the set indicator function that reads in the case of \eqref{dual} and \eqref{sub} as
\begin{equation}\label{nor-sub}
N_\O^-(\ox)=\sub^{-}\dd_{\O}(\ox)\;\mbox{ whenever }\;\ox\in\O.
\end{equation}
It follows from the above discussions that we always have the inclusions
\begin{equation}\label{sub-inc}
\Hat\partial\ph(\ox)\subset\partial^-\ph(\ox)\subset\Bar\partial\ph(\ox)
\end{equation}
between \eqref{sub} and the (convex) regular/Fr\'echet and convexified/Clarke subdifferentials, where the first inclusion (but not the second one) holds as equality in finite dimensions. Thus
\begin{equation}\label{lim-sub}
\partial^L \ph(\ox)=\Limsup_{x\st{\ph}{\to}\ox}\partial^{-} \ph(x)\;\mbox{ for all }\;\ox\in\dom\ph
\end{equation}
for the (nonconvex) limiting/Mordukhovich subdifferential of $\ph$ in $\X=\R^n$, where the symbol $x\st{\ph}{\to}\ox$ means that $x\to\ox$ with $\ph(x)\to\ph(\ox)$. As well known, the first inclusion in \eqref{sub-inc} and the representation \eqref{lim-sub} fail in infinite dimensions; see, e.g., $\ph(x):=-\|x\|$ with $\ox=0$ in the standard space ${\cal C}[0,1]$ of continuous functions, where $\Hat\partial\ph(\ox)=\partial^L\ph(\ox)=\emp$ while $\partial^-\ph(0)\ne\emp$.

Observe further that, by taking into account that $T_{\ss\epi\ph}(\ox,\ph(\ox))=\epi\d\ph(\ox)$, there exists an equivalent representation of this subdifferential via the (Dini-Hadamard) subderivative \eqref{subder}:
\begin{equation}\label{dhsub}
\sub^-\ph(\ox)=\big\{v\in\X^*\;\big|\;\la v,u\ra\le\d\ph(\ox)(u)\;\mbox{ for all }\;u\in X\big\},\quad\ox\in\dom\ph,
\end{equation}
which can be used as the subdifferential definition in this case.\vspace*{0.03in}

Finally in this section, we recall the two well-posedness properties of set-valued mappings, which have been well recognized in variational analysis and optimization. Given $F\colon\X\tto\Y$ with
\begin{equation*}
\dom F:=\big\{x\in\X\;\big|\;F(x)\ne\emp\big\}\;\mbox{ and }\;\gph F:=\big\{(x,y)\in\X\times\Y\;\big|\;y\in F(x)\big\},
\end{equation*}
it is said that $F$ is {\em metrically regular} around $(\ox,\oy)\in\gph F$ if there exist $\kappa>0$ and neighborhoods $U$ of $\ox$ and $V$ of $\oy$ such that the distance estimate
\begin{equation}\label{metreq}
\dist\big(x;F^{-1}(y)\big)\le\kappa\,\dist\big(y;F(x)\big)\;\mbox{ for all }\;(x,y)\in U\times V
\end{equation}
fulfills. If $y=\oy$ in \eqref{metreq}, the mapping $F$ is said to be {\em metrically subregular} at $(\ox,\oy)$. These and related properties have been broadly studied in variational analysis with numerous applications to optimization; see, e.g., the books \cite{bz,dr,i,m06,m18,penot,rw} and the references therein.\vspace*{-0.15in}

\section{Qualification Conditions}\label{sec03}\sce\vspace*{-0.05in}

In this section we define and investigate the qualification conditions used in the paper for the study and applications of the compositions given by
\begin{equation}\label{CS}
\ph(x):=\big(\theta\circ f\big)(x),\quad x\in\X,
\end{equation}
where $f\colon\X\to\Y$ is a single-valued mapping between normed spaces, and where $\theta\colon\Y\to\oR$. The next proposition formulates these qualification conditions and establishes relationships between them.\vspace*{-0.05in}

\begin{Proposition}[\bf relationships between qualification conditions]\label{equi} Dealing with composition \eqref{CS} in normed spaces, assume that $f\colon\X\to\Y$ is continuously differentiable around $\ox\in\X$, and that $\theta\colon\Y\to\oR$ is a continuous relative to its domain with  $\theta(f(\ox))\in\dom\th$. Consider the following assertions:\vspace*{-0.1in}
\begin{itemize}[noitemsep]
\item[\bf{(i)}] The mapping $H\colon\X\times\R\tto\Y\times\R$ defined by $H(x,\alpha):=(f(x),\alpha)-\epi\theta$ is metrically regular at $\big((\ox,\theta(f(\ox))),(0,0)\big)$.
\item[\bf{(ii)}] The mapping $H$ from {\rm(i)} is metrically subregular at this point.
\item[\bf{(iii)}] The mapping $G\colon\X\tto\Y$ defined by $G(x):=f(x)-\dom\theta$ is metrically subregular at $(\bar x,0)$.
\item[\bf{(iv)}] The {\sc Abadie qualification condition} {\rm(AQC)} holds at $\ox$, i.e.,
\begin{equation}\label{AQC}
T_{\ss\dom\ph}(\ox)=\big\{u\in\X\;\big|\;\nabla f(\ox)u\in T_{\ss\dom\theta}\big(f(\ox)\big)\big\}.
\end{equation}
\item[\bf{(v)}] The {\sc limiting Guignard qualification condition} holds at $\ox$, i.e.,
\begin{equation}\label{GuC}
N_{\ss\dom\ph}^-(\ox)=\overline{{\nabla f(\ox)^*N^{-}_{\ss\dom\theta}\big(f(\ox)\big)}}^*,
\end{equation}
where the bar$^*$ signifies the set closure in the weak$^*$ topology of $\Y^*$.
\end{itemize}\vspace*{-0.1in}
Then we always get that {\rm(i)}$\implies${\rm(ii)}$\implies ${\rm(iii)}$\implies${\rm(iv)}. If in addition the domain set $\dom\theta$ is convex, then we also have that {\rm(iv)}$\implies${\rm(v)}.
\end{Proposition}\vspace*{-0.2in}
\begin{proof} Implication (i)$\implies$(ii) is obvious. Although implication (ii)$\implies$(iii) was verified \cite[Proposition~3.1]{mms1} in finite dimension, a close look on the proof reveals that it works for any normed space. Applications (iii)$\implies$(iv) directly follows from the subderivative chain rule in Theorem~\ref{subdchain}(ii) with replacing therein the function $\theta$ by the indicator function $\delta_{\ss\dom\theta}$.

It remains to justify the last implication (iv)$\implies$(v) provided that $\dom\th$ is convex. Assuming (iv), pick $v\in N_{\ss\dom\theta}(f(\ox))$ and $u\in T_{\ss\dom\ph}(\ox)$. Then we get $\nabla f(\ox)u\in T_{\ss\dom\theta}(\ox)$, which implies in turn that
\begin{equation*}
\big\la\nabla f(\ox)^*v,u\big\ra=\big\la v,\nabla f(\ox)u\big\ra\le 0.
\end{equation*}
Since the latter holds for any $u\in T_{\ss\dom \ph}(\ox)$, it yields $\nabla f(\ox)^* v\in N_{\ss\dom\ph}^- (\ox)$. This tells us that
\begin{equation*}
\nabla f(\ox)^*N_{\ss\dom\theta}\big(f(\ox)\big)\subset N_{\ss\dom\ph}^-(\ox)\Longrightarrow\overline{{\nabla f(\ox)^*N^{-}_{\ss\dom\theta}\big(f(\ox)\big)}}^* \subset N_{\ss\dom\ph}^-(\ox).
\end{equation*}

To verify the opposite inclusion in \eqref{GuC}, pick any $v\in N_{\ss\dom\ph}^-(\ox)$ and then get by \eqref{dual} that $\la v,u\ra\le 0$ for all $u\in T_{\ss\dom\ph}(\ox)$. Hence we deduce from (v) that
\begin{equation}\label{equi1}
\la v,u\ra\le 0\;\mbox{ for all }\;u\in\X\;\mbox{ with }\;\nabla f(\ox)u\in T_{\ss\dom\theta}\big(f(\ox)\big).
\end{equation}
To show that $v\in\overline{{\nabla f(\ox)^*N^{-}_{\ss\dom\theta}\big(f(\ox)\big)}}^*$, assume the contrary and observe that the set $\overline{{\nabla f(\ox)^*N^{-}_{\ss\dom\theta}\big(f(\ox)\big)}}^*$ is weak$^*$ closed and convex in $\X^*$. Applying the strict convex separation theorem to the latter set and $\{v\}$ in the weak$^*$ topology of $\X^*$ gives us $\xi\in\X\setminus\{0\}$ and $\varepsilon>0$ such that
\begin{equation*}
\big\la\nabla f(\ox)\xi,y\big\ra=\big\la\xi,\nabla f(\ox)^*y\big\ra\le\la\xi,v\ra-\varepsilon\;\mbox{ for all }\;y\in N_{\ss\dom\theta}\big(f(\ox)\big).
\end{equation*}
Since $N_{\ss\dom\theta}(f(\ox))$ is a cone, it implies that $0\le\la\xi,v\ra-\varepsilon$ and hence
\begin{equation*}
\big\la\nabla f(\ox)\xi,y\big\ra\le 0\;\mbox{ whenever }\;y\in N_{\ss\dom\theta}\big(f(\ox)\big).
\end{equation*}
Using now the {\em full duality} between the tangent and normal cones to convex sets, we conclude that $\nabla f(\ox)\xi\in T_{\ss\dom\theta}(f(\ox))$. Thus it follows from \eqref{equi1} that $\la v,\xi\ra\le 0$, which contradicts the above assertion $\varepsilon\le\la v,\xi\ra$ and completes the proof of the proposition.
\end{proof}\vspace*{-0.05in}

It is well known (see, e.g., \cite[Theorem~2.84 and Proposition~2.89]{bs}) that, in the case where both spaces $\X$ and $\Y$ are Banach, the metric regularity of the mapping $G$ in Propositions~\ref{equi}(iii) is equivalent to the fulfillment of the {\em Robinson constraint qualification}
\begin{equation}\label{rob}
0\in\mbox{int}\,\big\{f(\ox)+\nabla f(\ox)\X-\dom\theta\big\},
\end{equation}
which is broadly used in optimization. Hence the metric subregularity qualification condition from Propositions~\ref{equi}(iii) is much weaker than the conventional one \eqref{rob}.\vspace*{0.03in}

The next example describes a finite-dimensional setting, where the Abadie qualification condition \eqref{AQC} is {\em strictly weaker} than the metric subregularity one from Proposition~\ref{equi}(iii) for compositions \eqref{CS} with closed and convex domain sets $\dom\th$. It shows furthermore that the closure operation in \eqref{GuC} cannot be generally removed even in finite dimensions.\vspace*{-0.05in}

\begin{Example}[\bf specifying relationships between qualification conditions in finite dimensions]\label{ACexample} Let $\Theta\subset\R^m$ be closed convex cone, let $A$ be a $n\times m$ matrix such that the set $A\Theta$ is not closed in $\R^n$, and let $\theta:=\delta_{\Theta^*}$ and $f(x):=A^*x$ in \eqref{CS}, where the $^*$-symbol indicates for the dual/polar cone for $\Th$ and the adjoint/transpose matrix for $A$. Such a pair $(A,\Th)$ exists and yields the followings properties:\vspace*{-0.1in}
\begin{itemize}[noitemsep]
\item[\bf{(i)}] The Abadie qualification condition \eqref{AQC} holds at $\ox=0$.
\item[\bf{(ii)}] The set $\nabla f(\ox)^*N_{\ss\dom\theta}(f(\ox))$ is not closed.
\item[\bf{(iii)}] The set-valued mapping $x\mapsto f(x)-\dom\theta$ is not metrically subregular at $\ox=0$.
\end{itemize}
\end{Example}\vspace*{-0.15in}
\begin{proof} First we verify all the conclusions in this example under the imposed assumptions and then construct a specific pair $(A,\Th)$ satisfying these assumptions.

Observe that $\nabla f(0)=A^*$, and that $T_{\Theta^*}(f(0))=T_{\Theta^*}(0)=\Theta$. To verify (i), define the set $\O:=\dom\ph=\{x\in\R^n\;|\;A^*x\in\Theta^*\}$, which is clearly a closed and convex cone. Thus we have
\begin{equation*}
T_{\O}(0)=\O=\big\{x\in\R^n\;\big|\;\nabla f(0)x\in T_{\Theta^*}(0)\big\}
\end{equation*}
while showing that AQC \eqref{AQC} is satisfied at $\ox=0$ for the composition under consideration.

To check further that the set in (ii) is not closed, observe that $N_{\Theta^*}(0)=\Theta$, and hence
\begin{equation}\label{aba}
\nabla f(0)^*N_{\ss\dom\theta}\big(f(0)\big)=AN_{\Theta^*}(0)=A\Theta,
\end{equation}
which verifies (ii) since the set $A\Th$ is assumed to be nonclosed in $\R^n$.

Turning finally to (iii), recall that the metric subregularity of the mapping $x\mapsto f(x)-\Theta^{*}$ at $\ox=0$ ensures the normal cone calculus rule
\begin{eqnarray}\label{aba1}
N_{\O}(0)=\nabla f(0)^*N_{\ss\dom\theta}\big(f(0)\big),
\end{eqnarray}
which is obtained in variational analysis even in more general settings; see, e.g., \cite[Lemma~2.1]{gm} and Theorem \ref{cahinset} below for the proofs and discussions. Thus the set on the right-hand side of \eqref{aba1} is closed that contradicts \eqref{aba} by the assumption on $A\Th$. This verifies (iii).\vspace*{0.03in}

To complete our consideration in this example, it remains to construct a closed convex cone $\Theta$ and a matrix $A$ satisfying the imposed assumptions. We do it in two stages. First note that there exist closed convex cones $\Theta_1,\Theta_2\subset\R^3$ such that their sum $\Theta_1+\Theta_2$ is not closed. Indeed, take
\begin{equation*}
\Theta_1:=\big\{(x,r)\in\R^2\times\R\;\big|\;\|x\|\le r\big\}\;\mbox{ and }\;\Theta_2:=\big\{t(1,0,-1)\in\R^3\;\big|\;t\ge 0\big\}.
\end{equation*}
The sum of these cones is not closed, since it does not contain $(0,1,0)$ while containing
\begin{equation*}
\overbrace{(-t,1+t^{-1},t+\frac{t^{-1}+2t^{-2}+t^{-3}}2)}^{\in\Theta_1}+\overbrace{(t,0,-t)}^{\in\Theta_2}=\Big(0,1+t^{-1},\frac{t^{-1}+2t^{-2}+t^{-3}}2\Big)
\end{equation*}
for all $t>0$. Define now the linear mapping $T\colon\R^3\times\R^3\to\R^3$ by $T(x,y):=x+y$ and observe that the set $\Theta_1\times\Theta_2$ is closed in $\R^3\times \R^3 $, but its image $T(\Theta_1\times\Theta_2)=\Theta_1+\Theta_2$ is not closed. This clearly gives us the pair $(A,\Th)$ satisfying the assumptions made.\vspace*{-0.05in}
\end{proof}

It is shown in Section~\ref{sec05a} that the replacement of the Abadie qualification condition \eqref{AQC} by the stronger metric subregularity one excludes the situation like in Example~\ref{ACexample} and brings us to the normal cone chain rule without the closure operation as in \eqref{GuC} even in fairly general infinite-dimensional settings. This is the basis for deriving {\em pointwise} optimal conditions of the KKT type in contrast to {\em asymptotic} ones generated by chain rules with closures.\vspace*{0.03in}

Finally, we formalize the major qualification condition, which is a {\em constraint qualification} in optimization problems, that is broadly used in this paper.\vspace*{-0.06in}

\begin{Definition}[\bf metric subregularity qualification condition]\label{defmscq} Let $\ph:=\th\circ f$ be a composition generated by a mapping $f\colon\X\to\Y$ between normed spaces and an extended-real-valued function $\th\colon\Y\to\oR$. We say that the {\sc metric subregularity qualification condition} $($MSQC$)$ holds for $\ph$ at $\ox\in\dom\ph$ with modulus $\kappa>0$ if the mapping $x\mapsto f(x)-\dom\theta$ is metrically subregular at $(\bar x,0)$ with this constant.
\end{Definition}\vspace*{-0.05in}

Using the notion of metric subregularity in \eqref{metreq} with $y=\oy=0$ therein, we reformulate Definition~\ref{defmscq} as follows: there exist a neighborhood $U$ of $\ox$  and a constant $\kappa>0$ such that the distance estimate
\begin{equation}\label{mscq}
\dist(x;\O)\le\kappa\,\dist\big(f(x);\dom\theta\big)\;\mbox{ with }\;\O:=\big\{x\in \X \;\big|\;f(x)\in\dom\theta\big\}
\end{equation}
holds for all $x\in U$. Note that MSQC \eqref{mscq} has been used in variational analysis and constrained optimization, mainly in more recent publications in finite-dimensional settings; see, e.g., \cite{g,gm,go,ho,io,mms1} and the references therein. Here we systematically employ it in infinite dimensions.
\vspace*{-0.15in}

\section{Calculus of Subderivatives in Normed Spaces}\sce\label{sec04}\vspace*{-0.05in}

This section is devoted to deriving chain and sum rules of the {\em equality type} for subderivatives \eqref{subder} in normed spaces. The main result is a chain rule for compositions \eqref{CS} obtained in the two general settings: under the MSQC from Definition~\ref{defmscq} and under AQC \eqref{AQC} combined with the epi-differentiability of the outer function in \eqref{CS}. We start with the following characterization of proper epi-differentiability for arbitrary extended-real-valued functions that is of its own interest.\vspace*{-0.05in}

\begin{Lemma}[\bf characterization of proper epi-differentiability]\label{epilemma} Let $\ph\colon\X\to\oR$ with $\ox\in\dom\ph$, and let $\psi\colon\X\to\oR$ be a proper extended-real-valued function $\psi\colon\X\to\oR$ such that $\psi(u)\le\d\ph(\ox)(u)$ for all $u\in\X$. Then $\ph$ is properly epi-differentiable at $\ox$ with $\psi(\cdot)=\d\ph(\ox)(\cdot)$ if and only if for every $u\in\X$ there exists a path $\Tilde u\colon[0,\varepsilon]\to\X$ $($not necessarily continuous$)$ with the properties $\lim_{t\dn 0}\Tilde u(t)=\Tilde u(0)=u$ and
\begin{equation}\label{ut}
\psi(u)=\lim_{t\dn 0}\frac{\ph\big(\ox+t\Tilde u(t)\big)-\ph(\ox)}{t},\quad u\in\X.
\end{equation}
\end{Lemma}\vspace*{-0.15in}
\begin{proof} Assume first that $\ph$ is epi-differentiable at $\ox$ with the subderivative function $\psi(\cdot)=\d\ph(\ox)(\cdot)$. If $u\in\X$ is such that $\psi(u)=\infty$, set $\Tilde u(t):=u$ for all $t\ge 0$ and observe that \eqref{ut} obviously holds. Suppose further that $\psi(u)$ is a finite number and then get $(u,\psi(u))\in T_{\ss\epi\ph}(\ox,\ph(\ox))$ by the definition of epi-differentiability. Moreover, the latter ensures that the pair $(u,\psi(u))$ is a derivable tangent in $T_{\ss\epi\ph}(\ox,\ph(\ox))$. This tells us that there exists a path $\xi\colon[0,\ve]\to\epi\ph$ with the components $\xi(t)=(\xi_1(t),\xi_2(t))$ for all $t\in[0,\ve]$ satisfying the conditions $\xi(0)=(\ox,\ph(\ox))$ and $\xi'_+(0)=(u,\psi(u))$. Setting now $\Tilde u(t):=\frac{\xi_1(t)-\ox}{t}$ for all $t\in(0,\ve]$ with $\Tilde u(0):=u$, we deduce from the inclusion $\xi(t)\in\epi\ph$ on $[0,\ve]$ that
\begin{equation*}
\frac{\ph\big(\ox+t\Tilde u(t)\big)-\ph(\ox)}{t}\le\frac{\xi_2(t)-\ph(\ox)}{t}\;\mbox{ for all }\;t\in(0,\ve].
\end{equation*}
This clearly implies that $\frac{\xi_2(t)-\ph(\ox)}{t}\to\psi(u)=\d\ph(\ox)(u)$ as $t\dn 0$ and shows that the limit of the quotient on the left-hand side above exists and is equal to $\psi(u)$.

To verify the opposite implication, suppose that \eqref{ut} holds and observe that it yields $\psi(\cdot)\ge\d\ph(\ox)(\cdot)$. Invoking the imposed assumption on $\psi(\cdot)\le\d\ph(\ox)(\cdot)$, this tells us that $\psi(\cdot)=\d\ph(\ox)(\cdot)$ and hence $\epi\psi=T_{\ss\epi\ph}(\ox,\ph(\ox))$. We need to show that every vector in $T_{\ss\epi\ph}(\ox,\ph(\ox))$ is a derivable tangent. To this end, pick any pair $(u,\alpha)\in T_{\ss\epi\ph}(\ox,\ph(\ox))=\epi\psi$ and let $\Tilde u\colon[0,\ve]\to\X$ be the path taken from \eqref{ut}. By setting finally $\xi\colon[0,\ve]\to\X \times \R$ with
\begin{equation*}
\xi(t):=\big(\ox+t\Tilde u(t),\ph(\ox+t \Tilde u(t))+t(\alpha-\psi(u))\big),\quad t\in[0,\ve],
\end{equation*}
it is easy to check that $\xi(t)\in\epi\ph$ for all $t\in[0,\ve]$, $\xi(0)=(\ox,\ph(\ox))$, and $\xi'_+(0)=(u,\psi(u))$, which therefore completes the proof of the lemma.
\end{proof}\vspace*{-0.05in}

Now we consider a large class of extended-real-valued functions, which play an important role in what follows. Recall that $\ph\colon\X\to\oR$ is {\em relatively Lipschitz continuous} around $\ox\in\dom\ph$ with respect to a set $\O\subset\dom\ph$ if there exist a constant $\ell\ge 0$ and a neighborhood $U$ of $\ox$ with
\begin{equation}\label{liprel}
|\ph(x)-\ph(u)|\le\ell\|x-u\|\;\mbox{ for all }\;x,u\in\O\cap U.
\end{equation}
This class is significantly broader than the standard class of locally Lipschitzian functions. As a trivial example, it include indicator functions of sets that plays a prominent role in constrained optimization. Another important subclass of relatively Lipschitzian functions is formed by {\em piecewise linear-quadratic} ones $\ph\colon\X\to\oR$ defined as follows: $\ph$ is such that $\dom\ph=\cup_{i=1}^{s}\O_i$, where $\O_i$ are polyhedral convex sets for $i=1,\ldots,s$, and $\ph$ admits the representation
\begin{equation}\label{PWLQ}
\ph(x)=B_i(x,x)+\langle b_i,x\rangle+\beta_i\;\mbox{ whenever }\;x\in\O_i,
\end{equation}
where $B_i\colon\X\times\X\to\R$ is a continuous and symmetric bilinear mapping, $b_i\in\X^*$, and $\beta_i\in\R$ for all $i=1,\ldots,s$. The next simple albeit useful lemma reveals a remarkable subderivative property of relatively Lipschitzian functions on normed spaces.\vspace*{-0.05in}

\begin{Lemma}[\bf subderivatives of relatively Lipschitzian functions]\label{subdomain} Let $\ph\colon\X\to\oR$ be relatively Lipschitzian around $\ox\in\dom\ph$ with respect to its domain. Then the subderivative function $\d\ph(\ox)\colon\X\to\oR $ is proper, and its domain is represented by $\dom\d\ph(\ox)=T_{\ss\dom\ph}(\ox)$.
\end{Lemma}\vspace*{-0.15in}
\begin{proof} Since the inclusion $\dom\d\ph(\ox)\subset T_{\ss\dom\ph}(\ox)$ always holds, let us verify the opposite one. Pick any $u\in T_{\ss\dom\ph}(\ox)$ and deduce from the imposed Lipschitzian property that for all $t>0$ and $u'\in\X$ with $\ox+t u'\in\dom\ph$ we have the estimate
\begin{equation*}
\Big|\frac{\ph(\ox+t u')-\ph(\ox)}{t}\Big|\le\ell\|u '\|,
\end{equation*}
where $\ell\ge 0$ is taken from \eqref{liprel}. This implies by the subderivative definition \eqref{subder} that
\begin{equation*}
|\d\ph(\ox)(u)|\le\ell\|u\|\;\mbox{ for all }\;u\in T_{\ss\dom\ph}(\ox),
\end{equation*}
which clearly shows that $\d\ph(\ox)$ is proper with $\d\ph(\ox)(0)=0$ and that $u\in\dom\d\ph(\ox)$.
\end{proof}\vspace*{-0.05in}

Now we are ready to establish the main result of this section, which provides the subderivative chain rule as equality under the aforementioned qualification conditions. We are not familiar with such results in infinite-dimensional spaces, while the finite-dimensional version of the chain rule in (ii) under the metric subregularity qualification condition was previously obtained in our paper with Sarabi \cite[Theorem~3.4]{mms1}. Assertion (i) of the following theorem is new even in finite dimensions.\vspace*{-0.05in}

\begin{Theorem}[\bf subderivative chain rule under the weakest qualification conditions]\label{subdchain} Let $f\colon\X\to\Y$ be a mapping between normed spaces, which is Fr\'echet differentiable at some point $\ox\in\X$, and let $\theta\colon\Y\to\oR$ be relatively Lipschitz continuous around $\oy:=f(\ox)$ with respect to its domain. Considering the composition $\ph=\theta \circ f$ in \eqref{CS}, assume that:\vspace*{-0.05in}
\begin{itemize}
\item[\bf{(i)}] either $\theta$ is epi-differentiable at $\oy$ and the Abadie qualification condition \eqref{AQC} holds at $\ox$,\vspace*{-0.1in}
\item[\bf{(ii)}] or the metric subregularity qualification condition \eqref{mscq} is satisfied at $\ox$.
\end{itemize}\vspace*{-0.05in}
Then we have the following subderivative chain rule for the composition in \eqref{CS}:
\begin{equation}\label{dchain}
\d(\theta\circ f)(\ox)(u)=\d\theta(\oy)\big(\nabla f(\ox)u\big)\;\mbox{ for all }\;u\in\X.
\end{equation}
Furthermore, the epi-differentiability of the function $\theta$ at $\oy$  in case {\rm(ii)} ensures that the composition $\ph=\theta\circ f$ is epi-differentiable at $\ox$.
\end{Theorem}\vspace*{-0.15in}
\begin{proof} First we show that the inequality ``$\ge$" in \eqref{dchain} holds whenever $f$ is Fr\'echet differentiable at $\ox$ without any other assumptions of the theorem. To proceed,
pick any $u\in\X$ and observe that $\nabla f(\ox)u'+\frac{o(t\|u'\|)}{t}\to\nabla f(\ox)u$ as $t\dn 0$ and $u'\to u$. This yields the relationships
\begin{eqnarray*}
\disp\d(\theta\circ f)(\ox)(u)&=&\liminf_{\substack{t\dn 0\\u'\to u}}\frac{\theta\big(f(\ox+t u')\big)-\theta(\oy)}{t}\nonumber\\
&=&\liminf_{\substack{t\dn 0\\u'\to u}}\frac{\theta\big(f(\oy+t\nabla f(\ox)u'+o(t\|w'\|)\big)-\theta(\oy)}{t}\nonumber\\
&=&\liminf_{\substack{t\dn 0\\u'\to u}}\frac{\theta\big(\oy+t\big(\nabla f(\ox)u'+\frac{o(t\|w'\|)}{t})\big)-\theta(\oy)}{t}\nonumber\\
&\ge&\d\theta(\oy)\big(\nabla f(\ox)u\big)\;\mbox{ for all }\;u\in\X,
\end{eqnarray*}
which verify the claimed inequality. To prove the opposite inequality in \eqref{dchain}, deduce from the Lemma~\ref{subdomain} that $\d\theta(\oy)(\nabla f(\ox)u)>-\infty$ for all $u\in\X$. Since the fulfillment \eqref{dchain} is obvious if $\d\theta(f(\oy)(\nabla f(\ox)u)=\infty$, we suppose in what follows that $\d\theta(\oy)(\nabla f(\ox)u)$ is finite. In case (i), where $\theta$ is properly epi-differentiable at $\oy$, Lemma~\eqref{epilemma} gives us a path $y(\cdot)$ in $\Y$ such that $\lim_{t\dn 0}y(t)=\nabla f(\ox)u$ and
\begin{equation}\label{subd0}
\d\theta(\oy)\big(\nabla f(\ox)u\big)=\lim_{\substack{t\dn 0}}\frac{\theta\big(\oy+ty(t)\big)-\theta(\oy)}{t}.
\end{equation}
Since $\d\theta(\oy)(\nabla f(\ox)u)$ is finite, we suppose without lost of generality that $\oy+t y(t)\in\dom\theta $ for all $t\in[0,\ve]$, and hence $\nabla f(\ox)u\in\dom\d\theta(\oy)=T_{\ss\dom\theta}(\oy)$ by Lemma~\ref{subdomain}. Thus the assumed AQC \eqref{AQC} yields $u\in T_{\ss\dom\ph}(\ox)$. This gives us by definition \eqref{tan} sequences
$\{u_{k}\}\subset\X$ and $t_k>0$ converging to $u$ and $0$, respectively, such that $\ox+t_k u_k\in\dom\ph$ for all $k\in\N$. Using these along with \eqref{subd0} leads us to
\begin{eqnarray*}
\d\theta(\oy)\big(\nabla f(\ox)u\big)&=&\disp\lim_{k\to\infty}\bigg[\frac{\theta\big(f(\ox+t_ku_k)\big)-\theta(\oy)}{t_k}+\frac{\theta\big(\oy+t_k y(t_k)\big)-\theta\big(f(\ox+t_ku_k)\big)}{t_k}\bigg] \nonumber\\\nonumber
&\ge&\disp\liminf_{k\to\infty}\frac{\theta\big(f(\ox+t_ku_k)\big)-\theta(\oy)}{t_k}-\ell\limsup_{k\to\infty}\Big\|\frac{f(\ox+t_ku_k)-\oy}{t_k}-y(t_k)\Big\|\\
&\ge&\d(\theta\circ f)(\ox)(u)-\ell\limsup_{k\to\infty}\Big\|\nabla f(\ox)u_k +\frac{o(t_k)}{t_k}-y(t_k)\Big\|,
\end{eqnarray*}
where $\ell$ is a constant of the relative Lipschitz continuity of $\theta$ around $\oy$ with respect to its domain. Since $y(t_k)\to\nabla f(\ox)u$ as $k\to\infty$, we get $\limsup_{k\to\infty}\|\nabla f(\ox)u_k+\frac{o(t_k)}{t_k}-y(t_k)\|=0$, which verifies the subderivative chain rule \eqref{dchain} under the assumptions in (i).

The subderivative chain rule \eqref{dchain} under \eqref{mscq} in (ii) was proved in \cite[Theorem~3.4]{mms1} in finite-dimensional spaces, but the proof therein works in the case of normed spaces as well.

To complete the proof of the theorem, it remains to show that the epi-differentiability $\theta$ at $\oy$ in case (ii) yields this property of the composition $\ph$ at $\ox$, which is new even in finite dimensions. To proceed, deduce first from Lemma~\ref{subdomain} that the function $\psi(u):=\d\theta(\oy)(\nabla f(\ox)u)=\d\ph(\ox)(u)$ is proper. Then applying Lemma~\ref{epilemma} verifies the conclusion of the theorem provided the fulfillment of representation \eqref{ut} for the function $\psi$ whenever $u\in\X$. To this end, pick $u\in\X$ with $\d\psi (u) = \theta(\oy)(\nabla f(\ox)u) < \infty $ and get by Lemma~\ref{subdomain} that $\nabla f(\ox)u\in T_{\ss\dom\theta}$. Choose further a path $y(\cdot)$ in $\Y$ for which \eqref{subd0} holds. Then applying MSQC \eqref{mscq} with some $\kappa > 0$ and $x:=\ox+t u$ for sufficiently small $t>0$ gives us
\begin{equation*}
\dist\big(\ox+t u;\dom\ph\big)\le\kappa\,\dist\big(f(\ox+t u);\dom\th\big),
\end{equation*}
which in turn results in the estimates
\begin{eqnarray*}
\dist\Big(u;\frac{\dom\ph-\ox}{t}\Big)&\le&\frac{\kappa}{t}\,\dist\big(\oy+t\nabla f(\ox)u+o(t);\dom\theta\big)\nonumber\\
&\le&\frac{\kappa}{t}\,\Big\|\oy+t\nabla f(\ox)u+o(t)-\oy-t y(t)\Big\|\nonumber\\
&=&\kappa\,\Big\|\nabla f(\ox)u-y(t)+\frac{o(t)}{t}\Big\|.
\end{eqnarray*}
Thus for all $t>0$ sufficiently small we find $\Tilde u(t)\in\frac{\dom\ph-\ox}{t}$ such that
\begin{equation*}
\|u-\Tilde u(t)\|\le\kappa\,\Big\|\nabla f(\ox)u-y(t)+\frac{o(t)}{t}\Big\|+t.
\end{equation*}
This tells us that $\ox+t\Tilde u(t)\in\O$ for all small $t>0$, and that $\lim_{t\dn 0}\Tilde u(t)=u$. Since $\theta$ is relatively Lipschitz continuous around $\oy$ with respect of its domain, we get
\begin{equation*}
\Big\|\frac{\theta\big(\oy+t y(t)\big)-\theta\big(f(\ox+t\Tilde u(t))\big)}{t}\Big\|=\ell\,\Big\|\nabla f(\ox)\Tilde u(t)+\frac{o(t)}{t}-y(t)\Big\|\to 0 \;\mbox{ as }\;t\dn 0.
\end{equation*}
Taking all the above into consideration, we arrive at the equalities
\begin{eqnarray*}
\psi(u)&=&\disp\lim_{t\dn 0}\bigg[\frac{\theta\big(f(\ox+t\Tilde u(t))\big)-\theta(\oy)}{t}+\frac{\theta\big(\oy+t y(t)\big)-\theta\big(f(\ox+t\Tilde u(t))\big)}{t}\bigg]\nonumber\\\nonumber
&=& \disp\lim_{t\dn 0}\frac{\theta\big(f(\ox+t\Tilde u(t))\big)-\theta(\oy)}{t}=\lim_{t\dn 0}\frac{\ph\big(\ox+t\Tilde u(t)\big)-\ph(\ox)}{t},
\end{eqnarray*}
which verify \eqref{ut} and thus complete the proof of the theorem.
\end{proof}\vspace*{-0.05in}

Finally in this section, we establish the subderivative {\em sum rule} for extended-real-valued functions on normed spaces under the fulfillment of the corresponding summation counterparts of the Abadie and metric subregularity qualification conditions.\vspace*{-0.05in}

\begin{Theorem}[\bf subderivative sum rule]\label{subdsum} Let $\phi\colon\X\to\oR$ and $\psi\colon\X\to\oR$ be defined on a normed space $\X$, and let $\ox\in(\dom\phi)\cap(\dom\psi)$. Assume that both $\phi$ and $\phi$ are relatively Lipschitzian around $\ox$ with respect to their domains and are epi-differentiable at this point. Suppose also that
either\vspace*{-0.1in}
\begin{enumerate}
\item[\bf{(i)}] the tangential qualification condition
\begin{equation}\label{tancq}
T_{\ss(\dom\phi)\cap(\dom\psi)}(\ox)=T_{\ss\dom\phi}(\ox)\cap T_{\ss\dom\psi}(\ox)
\end{equation}
is satisfied at $\ox$, or the following\vspace*{-0.1in}
\item[\bf{(ii)}] metric qualification condition holds at this point: there exist a number $\kappa>0$ and a neighborhood $U$ of $\ox$ ensuring the fulfillment of the distance estimate
\begin{equation}\label{mscqsum}
\dist\big(x;\dom\phi\cap\dom\psi\big)\le\kappa\,\big(\dist(x;\dom\phi)+\dist(x;\dom\psi)\big)\;\mbox{ for all }\;x\in U.
\end{equation}
\end{enumerate}\vspace*{-0.1in}
Then we have the subderivative sum rule
\begin{equation}\label{subder-sum}
\d(\phi+\psi)(\ox)(u)=\d\phi(\ox)(u)+\d\psi(\ox)(u),\quad u\in\X.
\end{equation}
Furthermore, in case {\rm(ii)} the sum $\phi+\psi$ is epi-differentiable at $\ox$.
\end{Theorem}\vspace*{-0.15in}
\begin{proof}
Setting $\Y:=\X\times\X$, define $f\colon\X\to\Y$ by $f(x):=(x,x)$ and $\th\colon\Y\to\oR$ by $\th(x,z):=\phi(x)+\psi(z)$. We clearly have that $\phi+\psi=\theta\circ f$, and that $\theta$ is relatively Lipschitzian around $(\ox,\ox)$ with respect to the domain $\dom\phi\times\dom\psi$. Let us first show that $\theta$ is epi-differentiable at $(\ox,\ox)$ and the representation
\begin{equation}\label{th-sum}
\d\th(\ox,\ox)(u,z)=\d\phi(\ox)(u)+\d\psi(\ox)(z),\quad(u,z)\in\Y,
\end{equation}
holds even without imposing the qualification conditions from (i) and (ii). Indeed, it is obvious that
\begin{equation*}
\Delta_t\th(\ox,\ox)(u,z)=\Delta_t\phi(\ox)(u)+\Delta_t\psi(\ox)(z)\;\mbox{ whenever }\;(u ,z)\in\Y:=\X\times\X
\end{equation*}
for the finite difference $\Delta_t$ in \eqref{Delta}. Fixing $(u,z)\in\Y$ and applying the epi-differentiability criterion of Lemma~\ref{epilemma} to both functions $\phi$ and $\psi$, we find the corresponding paths $\tilde{u}(\cdot)$ and $\tilde{z}(\cdot)$ such that
\begin{eqnarray*}
\disp\d\theta(\ox,\ox)(u,z)&\le&\liminf_{t\dn 0}\Delta_t\theta(\ox,\ox)\big(\tilde{u}(t),\tilde{z}(t)\big)=\liminf_{t\dn 0}\big[\Delta_t\phi(\ox)\big(\tilde{u}(t)\big)+\Delta_t\psi(\ox)\big(\tilde{z}(t)\big)\big]\\
&=&\lim_{t\dn 0}\Delta_t\phi(\ox)\big(\tilde{u}(t)\big)+\lim_{t\dn 0}\Delta_t\psi(\ox)\big(\tilde{z}(t)\big)=\d\phi(\ox)(u)+\d\psi(\ox)(z).
\end{eqnarray*}
This proves the inequality ``$\le$" in \eqref{th-sum} and that $\th$ is is epi-differentiable at $(\ox,\ox)$ by Lemma~\ref{epilemma}. Moreover, we have the equality in \eqref{th-sum}, since the opposite inequality ``$\ge$" therein is trivial.

To proceed  further with the proof of \eqref{subder-sum}, it is not hard to check that the qualification conditions \eqref{tancq} and \eqref{mscqsum} reduces to AQC \eqref{AQC} and MSQC \eqref{mscq}, respectively, for the composition $\theta\circ f$. Applying now Theorem~\ref{subdchain}, for all $u\in\X$ we get the equalities
\begin{eqnarray*}
\d(\ph+\psi)(\ox)(u)&=&\d(\th\circ f)(\ox)(u)=\d\th(\oy)\big(\nabla f(\ox)u\big)\\
&=&\d\th(\ox,\ox)(u,u)=\d\phi(\ox)(u)+\d\psi(\ox)(u),
\end{eqnarray*}
which verifies the sum rule \eqref{subder-sum}. Since in case (ii) we have MSQC \eqref{mscq} for the composition $\th\circ f$ defined above, and since we have already proved that the epi-differentiability of both $\phi$ and $\psi$ at $\ox$ yields the epi-differentiability of $\th$ at $(\ox,\ox)$, we deduce the epi-differentiability of the sum $\phi+\psi$ at $\ox$ from Theorem~\ref{subdchain}. This completes the proof of the theorem.
\end{proof}
\vspace*{-0.25in}

\section{Subdifferential Calculus}\sce \label{sec05}\vspace*{-0.05in}

This section concerns some calculus rules for the Dini-Hadamard subdifferential \eqref{sub} of extended-real-valued functions defined on normed spaces. We obtain here new subdifferential chain and sum rules in the form of {\em equalities}. Besides being important for their own sake, such equality-type results play a crucial role for further developments in second-order variational analysis.

It has been realized in variational analysis and generalized differentiation that the equality-type calculus rules for various subdifferentials involve appropriate regularity conditions. A deep investigation of different regularity notions and relationships between them in generalized differential calculus was done by Bounkhel and Thibault; see \cite{b,bt} and the references therein. Now we formulate a natural notion of Dini-Hadamard regularity for extended-real-valued functions and set in normed spaces that is suitable for our subsequent study in this paper.\vspace*{-0.05in}

\begin{Definition}[\bf Dini-Hadamard regularity]\label{DHreg} Let $\ph\colon\X\to\oR$ be an extended-real-valued function on a normed space $\X$, and let $\ox\in\dom\ph$. We say that $\ph$ is {\sc Dini-Hadamard regular} at $\ox$ if
\begin{equation}\label{reg1}
\sup\big\{\la v,u\ra\;\big|\;v\in\sub^-\ph(\ox)\big\}=\d\ph(\ox)(u)\;\mbox{ for all }\;u\in\X.
\end{equation}
The set $\Omega\subset\X$ is Dini-Hadamard regular at $\ox\in\O$ if the indicator function
$\delta_\O$ has this property.
\end{Definition}\vspace*{-0.05in}

Note that if $\ph$ is either Fr\'echet differentiable at $\ox$, or locally Lipschitzian around $\ox$ and G\^ateaux differentiable at this point, then it is Dini-Hadamard regular at $\ox$. This is due to the classical relationship $\d\ph(\ox)(u)=\la\nabla\ph(\ox),u\ra$ for all $u\in\X$ under such a differentiability. Furthermore, we clearly have that any convex function and set is regular in the sense of Definition~\ref{DHreg}. In the general nonconvex case even in finite dimensions, the Dini-Hadamard regularity \eqref{reg1} is {\em strictly weaker} than the commonly used subdifferential regularity by Clarke \cite{c,rw} and lower regularity by Mordukhovich \cite{m06,m18}. In particular, the function $\ph\colon\R\to\R$ defined by $\ph(x):=x^{2}\sin\frac{1}{x}$ if $x\ne 0$ and $\ph(0):=0$ is Dini-Hadamard regular at $\ox=0$ while both Clarke and Mordukhovich regularities fail. It is worth mentioning that seemingly restrictive regularity conditions are rather appropriate for constrained optimization. In particular, they hold automatically for problems of conic programming that are reduced to minimizing the composite functions of type \eqref{CS}, where $\theta$ is convex and locally Lipschitz relative to $\dom\theta$, and where $f$ is a smooth function; cf.\ \cite{mms1,mms} and Section~\ref{sect06} below for more details.\vspace*{0.05in}

The next subdifferential chain rule is the main result of this section.\vspace*{-0.03in}

\begin{Theorem}[\bf Dini-Hadamard subdifferential chain rule]\label{lipchainrule} Let $f\colon\X\to\Y$ be a mapping between normed spaces that is continuously differentiable around $\ox$, and let $\theta\colon\Y\to\oR$ be locally Lipschitzian around $\oy:=f(\ox)$ and Dini-Hadamard regular at this point. Then the composition $\ph=\theta\circ f$ is Dini-Hadamard regular at $\ox$, and we have the following chain rule as equality:
\begin{equation}\label{lipchain}
\sub^-(\theta\circ f)(\ox)=\nabla f(\ox)^*\sub^-\theta(\oy).
\end{equation}
\end{Theorem}\vspace*{-0.15in}
\begin{proof} It follows from the smoothness of $f$ around $\ox$ that $\oy\in\mbox{int}(\dom\theta)$, and thus MSQC \eqref{mscq} holds at $\ox$ for the composition \eqref{CS}. This ensures by Theorem~\ref{subdchain}(ii) the fulfillment of the subderivative chain rule \eqref{dchain}. To verify now the inclusion ``$\supset$" in \eqref{lipchain}, pick any $v\in\sub^-\theta(\oy)$ and $u\in\X$ and then deduce directly from \eqref{dchain} and the definitions that
\begin{equation*}
\la\nabla f(\ox)^*v,u\ra=\la v,\nabla f(\ox)^*u\ra\le\d\theta(\oy)\big(\nabla f(\ox)u\big)=\d(\theta\circ f)(\ox)(u),
\end{equation*}
which yields $\nabla f(\ox)^*v\in\sub^-(\theta\circ f)(\ox))$ and thus justifies the claimed inclusion.

To verify the opposite inclusion ``$\subset$" in \eqref{lipchain}, suppose on the contrary that there exists $v\in\sub^-(\theta\circ f)(\ox)$ such that $v\notin\nabla f(\ox)^*\sub^-\theta(\oy)$. The assumed Lipschitz continuity of $\theta$ around $\oy$ ensures that the subgradient set $\sub^-\theta(\oy)$ is bounded and weak$^*$ closed in $\Y^*$, and thus it is compact in the weak$^*$ topology of $\Y^*$ by the Banach-Alaoglu theorem. Since $\nabla f(\ox)^*\colon\Y^*\to\X^*$ is linear and bounded operator, we conclude that the set $\nabla f(\ox)^*\sub^-\theta(\oy)$ is convex and weak$^*$ compact in $\X^*$. The strict convex separation theorem gives us $\xi\in\X\setminus\{0\}$ and $\varepsilon>0$ such that
\begin{equation*}
\la\xi,\nabla f(\ox)^*\lambda\ra\le\la\xi,v\ra-\varepsilon\;\mbox{ for all }\;\lambda\in\sub^-\theta(\oy).
\end{equation*}
Employing this together with the imposed Dini-Hadamard regularity of $\th$ at $\oy$ and the subderivative chain rule \eqref{dchain}, we arrive at the relationships
\begin{eqnarray*}
\disp\langle\xi,v\rangle-\varepsilon&\ge&\sup_{\lambda\in\sub^-\theta(\oy)}\langle\xi,\nabla f(\ox)^*\lambda\rangle\nonumber\\
&=&\sup_{\lambda\in\sub^-\theta(\oy)}\langle\nabla f(\ox)\xi,\lambda\rangle\nonumber\\
&=&\d\theta(\oy)\big(\nabla f(\ox)\xi\big)=\d(\theta\circ f)(\ox)(\xi)
\end{eqnarray*}
ensuring that $\langle\xi,v\rangle>\d(\theta\circ f)(\ox)(\xi)$. This tells us that $v\notin\sub^-(\theta\circ f)(\ox)$, a contradiction that completes the proof of the subdifferential chain rule \eqref{lipchain}.

It remains to verify that $\th\circ f$ is Dini-Hadamard regular at $\ox$ under the assumptions made. To proceed, take any $u\in\X$ and then deduce from Definition~\ref{DHreg} and the subderivative chain rule \eqref{dchain} that
\begin{eqnarray*}
\disp\sup_{v\in\sub^-(\theta\circ f)(\ox)}\langle v,u\rangle&=&\sup_{\lambda\in\sub^-\theta(\oy)}\langle\nabla f(\ox)^*\lambda,u\rangle=\sup_{\lambda\in\sub\theta(\oy)}\langle\lambda,\nabla f(\ox)u\rangle\\
&=&\d\theta(\oy)\big(\nabla f(\ox)u\big)=\d(\theta\circ f)(\ox)(u),
\end{eqnarray*}
which therefore finishes the proof of the theorem.
\end{proof}\vspace*{-0.05in}

It is not difficult to derive from Theorem~\ref{lipchainrule} the following subdifferential {\em sum rule} as {\em equality}.\vspace*{-0.05in}

\begin{Corollary}[\bf Dini-Hadamard subdifferential sum rule]\label{subdisum} Let $\X$ be a normed space, and let $\phi,\psi\colon\X\to\R$ be Lipschitz continuous around $\ox$. Assume that both $\phi$ and $\phi$ are Dini-Hadamard regular at $\ox$. Then the sum $\phi+\psi$ is also Dini-Hadamard regular at this point, and we have
\begin{equation*}\label{subdisum1}
\sub^-(\phi+\psi)(\ox)=\sub^-\phi(\ox)+\sub^-\psi(\ox).
\end{equation*}
\end{Corollary}\vspace*{-0.15in}
\begin{proof}
Define $\th\colon\X\times\X\to\R$ by $\th(x,y):=\phi(x)+\psi(y)$ and $f\colon\X\to\X\times\X$ by $f(x):=(x,x)$. It is easy to deduce directly from the subdifferential and regularity definitions that the Dini-Hadamard regularity of $\phi$ and $\psi$ at $\ox$ yields this property for the function $\th$ at $\oy:=(\ox,\ox)$, and furthermore we get
\begin{equation*}
\sub^-\th(\ox,\ox)=\sub^-\phi(\ox)\times\sub^-\psi(\ox).
\end{equation*}
Observe the representation $\phi+\psi=\th\circ f$. Applying now Theorem~\ref{lipchainrule} to this composition with $\oy:=(\ox,\ox)$ tells us that the sum $\phi+\psi$ is Dini-Hadamard regular at $\ox$ and that
\begin{eqnarray*}
\disp\sub^-(\phi+\psi)(\ox)=\sub^-(\th\circ f)(\ox)=\nabla f(\ox)^*\sub^-\th(\oy)=\sub^-\phi(\ox)+\sub^-\psi(\ox),
\end{eqnarray*}
which therefore completes the proof of the corollary.
\end{proof}\vspace*{-0.05in}

The next result establishes the subdifferential chain rule \eqref{lipchain} in the case of extended-real-valued and convex outer functions under the Abadie qualification condition.\vspace*{-0.05in}

\begin{Theorem}[\bf subdifferential chain rule under AQC]\label{2chain} Let $f\colon\X\to\Y$ be a mapping between normed spaces that is assumed to be Fr\'echet differentiable at $\ox\in\X$, and let $\theta\colon\Y\to\oR$ be convex and relatively Lipschitz continuous around $\oy=f(\ox)$ with respect to its domain. Suppose also that the Abadie qualification condition \eqref{AQC} holds at $\ox$, and that the set $\nabla f(\ox)^*\partial^- \theta(\oy)$ is weak$^*$ closed in $X^*$. Then we have the Dini-Hadamard subdifferential chain rule \eqref{lipchain}.
\end{Theorem}\vspace*{-0.15in}
\begin{proof} Observe that the usage of the full (vs.\ relative) Lipschitz continuity of $\th$ around $\oy$ allows us to verify in the proof of Theorem~\ref{lipchainrule} that the set $\nabla f(\ox)^*\partial^- \theta(\oy)$ is weak$^*$ closed in $\X^*$. Now we postulate the latter property and can employ assertion (i) of Theorem~\ref{subdchain}, which gives us the subderivative chain rule \eqref{dchain} under the Abadie qualification condition \eqref{AQC} and the epi-differentiability of $\th$ at $\oy$. Note that the imposed convexity of $\th$ ensures both Dini-Hadamard regularity and epi-differentiability of this function at $\oy$. The rest of the proof of the theorem follows the lines in the proof of Theorem \ref{lipchainrule}.
\end{proof}\vspace*{-0.05in}

As a consequence of Theorem~\ref{2chain} and the previous result of Theorem~\ref{subdchain} for subderivatives, we describe now a class of functions on normed spaces for which both subdifferential and subderivative chain rules hold unconditionally. The compositions considered below involve piecewise linear-quadratic functions \eqref{PWLQ} and frequently appear in constrained optimization. In finite dimensions, such compositions belong (under the additional convexity assumption for outer functions) to the family of {\em fully amenable functions} \cite{rw}, and the chain rules presented below can be derived by using MSQC \eqref{mscq} due to the classical Hoffman lemma; see \cite[Corollary~3.8]{mms1} and compare with the guides to \cite[Exercise~10.22(b)]{rw} that propose to employ different, more demanding results of finite-dimensional variational analysis. In the normed space setting under consideration we relay instead on the Abadie qualification condition.\vspace*{-0.05in}

\begin{Corollary}[\bf unconditional chain rules for remarkable compositions]\label{chpwlq} Let $\ph(x):=\theta(Ax+a)$, where $\theta\colon\Y\to\oR$ is a piecewise linear-quadratic function, $A\colon\X\to\Y$ is a linear bounded operator between normed space, and $a\in\Y$. Then we have
\begin{equation}\label{Achain1}
\d\ph(x)(w)=\d\theta(Ax+a)(Aw)\;\mbox{ for all }\;x\in\dom\ph\;\mbox{ and all }\;w\in\X.
\end{equation}
If in addition the outer function $\th$ is convex and the set $A^*\partial\th(Ax)$ is weak$^*$ closed in in $Y^*$, then
\begin{equation}\label{Achain2}
\sub^-\ph(x)=A^*\sub\theta(Ax+a)\;\mbox{ for all }\;x\in\dom\ph.
\end{equation}
\end{Corollary}\vspace*{-0.15in}
\begin{proof} As discussed in Section~\ref{sec04}, for any $x\in\dom\ph$ with $Ax\in\dom\th$ we have that the piecewise linear-quadratic function $\th$ is relatively Lipschitzian around $Ax\in\dom\th$ with respect to its domain. Furthermore, this function is clearly epi-differentiable at the point $y:=Ax$ and the Abadie qualification condition holds for $\ph$ at such $x$. Thus the subderivative chain rule \eqref{Achain1} is a consequence of Theorem~\ref{dchain}. The imposed additional assumptions allow us to deduce the subdifferential chain rule \eqref{Achain2} directly from Theorem~\ref{2chain}. Note that the closedness of $A^*\partial\th(Ax)$ is automatic if $\dim(\X)\times\dim(\Y)<\infty$ since the polyhedrality of the set $\partial\th(Ax)$ implies the polyhedrality (and hence closedness) of $A^*\partial\th(Ax)$.
\end{proof}\vspace*{-0.15in}

\begin{Remark}[\bf comparison with known results]\label{compar} {\rm If both spaces $\X$ and $\Y$ are finite-dimensional, Theorem~\ref{subdisum} reduces to the well-known results on subdifferential calculus under the imposed subdifferential regularity, which is the same for all the aforementioned subgradient mappings of locally Lipschitzian functions; see, e.g., \cite{m18,rw}. If both spaces $\X$ and $\Y$ are Banach and $f=A\colon \X \to \Y $ is a linear, bounded, and surjective operator, the chain rule \eqref{lipchain} is also well known; see, e.g., \cite[Proposition~4.29]{i}. If $\X$ is an (incomplete) normed space, the only version of \eqref{lipchain} we are familiar with can be found in Penot \cite[Proposition~4.42]{penot}, where $\dim(\Y)<\infty$ and $\nabla f(\ox)\colon\X\to\Y$ is surjective. From the viewpoint of applications to optimization, the above surjectivity assumption is quite restrictive, while the result of Theorem~\ref{lipchainrule} leads us to more reliable necessary optimality conditions; see Sections~\ref{sect06} and \ref{sect07}. The subderivative chain rule obtained in Theorem~\ref{2chain} is new even in finite dimensions.}
\end{Remark}\vspace*{-0.25in}

\section{Subamenable Functions and Sets in Infinite Dimensions}\sce\label{sec05a}\vspace*{-0.05in}

In this section we introduce and investigate a new class of extended-real-valued functions (and hence sets) in infinite dimensions, which plays a significant role in the subsequent material.\vspace*{-0.05in}

\begin{Definition}[\bf subamenable functions and sets]\label{subamenable} Let $\ph\colon\X\to\oR$ be a extended-real-valued function on a normed space $\X$, and let $\ox\in\dom\ph$. We say that $\ph$ is {\sc subamenable} at $\ox$ if there exist a neighborhood $U$ of $\ox$, a continuously differentiable mapping $f\colon\X\to\Y$ with values in a normed space, and a proper convex l.s.c.\ function $\theta\colon\Y\to\oR$, which is relatively Lipschitzian around $f(\ox)$ with respect to its domain $\dom\theta$, such that $\ph$ admits the representation
\begin{equation*}
\ph(x)=(\theta\circ f)(x)\;\mbox{ for all }\;x\in U
\end{equation*}
under the fulfillment of MSQC \eqref{mscq} at $\ox$. Correspondingly, a set $\O\subset\X$ is subamenable at $\ox\in\O$ if it is closed and its indicator $\delta_{\ss\O}$ is a subamenable function at this point.
\end{Definition}\vspace*{-0.07in}

In the case of sets in finite-dimensional spaces, this notion was introduced by Gfrerer and Mordukhovich in \cite[Definition~2.2]{gm} and applied there to Robinson stability of constraint systems. A fully subamenable version of this property for functions (where $\th\colon\R^m\to\oR$ is convex piecewise linear-quadratic, and where $f\colon\R^n\to\R^m$ is ${\cal C}^2$-smooth) was defined in our paper with Sarabi \cite[Definition~4.2]{mms1} with applications to finite-dimensional optimization. The origin of such properties goes back to Rockafellar \cite{r88} who introduced the {\em amenable} prototype of Definition~\ref{subamenable} and its modifications with $\th\colon\R^m\to\oR$ being lower semicontinuous and satisfying the {\em metric regularity} qualification condition (via its characterization) instead of MSQC \eqref{mscq}. Then broad applications of the amenability concepts to finite-dimensional variational analysis and optimization have been strongly developed in numerous publications; see \cite{rw} for more details and also \cite{m18} for recent updates and references.\vspace*{0.03in}

In this paper we proceed with the defined subamenability notion in infinite-dimensional spaces. It follows directly from Theorem~\ref{subdchain} that if $\ph\colon\X\to\oR$ on a normed space $\X$ is subamenable at $\ox\in\dom\ph$, then it is epi-differentiable at this point and satisfies the subderivative chain rule \eqref{dchain}. Note further that the subamenability of a set $\O$ at $\ox\in\O$ amounts to the existence of a neighborhood $U$ of $\ox$, a continuously differentiable mapping $f\colon\X\to\Y$, and a closed convex set $\Theta$ such that $\O$ admits the representation
\begin{equation}\label{Cset}
\O\cap U=\big\{x\in U\;\big|\;f(x)\in\Theta\big\}
\end{equation}
under the fulfillment of the metric subregularity qualification condition \eqref{mscq} at $\ox$ written now as
\begin{equation}\label{msqcset}
\dist(x;\O)\le\kappa\,\dist\big(f(x);\Theta\big)\;\mbox{ for all }\;x\in U
\end{equation}
with some constant $\kappa>0$. Then Theorem~\ref{subdchain} states that any set $\O\subset\X$ subamenable at $\ox$ is derivable at this point, and that its tangent cone at $\ox$ is calculated by
\begin{equation}\label{tan-msr}
T_{\O}(\ox)=\big\{u\in\X\;\big|\;\nabla f(\ox)u\in T_{\Theta}(\oy)\big\}\;\mbox{ with }\;\oy:=f(\ox).
\end{equation}
Furthermore, it is not hard to see that subamenability of $\O$ defined in \eqref{Cset} and \eqref{msqcset} is a {\em stable} property in the sense that if $\O$ is subamenable at $\ox$, then there exists a neighborhood $\Tilde U$ of $\ox$ such that $\O$ is subamenable at $x$ for all $x\in\O\cap\Tilde U$.\vspace*{0.03in}

The next important proposition reveals a remarkable property of the (sequential) {\em weak tangent cone} to a {\em subamenable} set in an arbitrary Banach space $\X$. It simply says that the weak sequential convergence in $\X$ does not chance the tangent cone \eqref{tan} to such sets. Note that this is no longer true if the sequential convergence is replaced by the topological/net one, which does not yields boundedness.\vspace*{-0.05in}

\begin{Proposition}[\bf tangent and weak tangent cones to subamenable sets]\label{weaktan} Let $\X$ be a normed space, and let $\O\subset\X$ be subamenable at $\ox\in\O$. Then we have
\begin{equation}\label{wtan}
T_\O (\ox)=T_{\O}^w(\ox):=\big\{u\in\X\;\big|\;\exists\,t_k{\downarrow}0,\;u_k\wto u\;\mbox{ as }\;k\to\infty\;\;\mbox{with}\;\;\ox+t_k u_k\in\O\big\}.
\end{equation}
\end{Proposition}\vspace*{-0.05in}
\begin{proof}
Taking into account that the inclusion $T_\O(\ox)\subset T_{\O}^w(\ox)$ is trivial for any set, we need verifying the opposite inclusion for subamenable ones. To proceed, use \eqref{Cset} and pick $u\in T_{\O}^w(\ox)$. It tells by the left-hand side of \eqref{wtan} that there exist sequences $t_k\dn 0$ and $u_k\wto u$ as $k\to\infty$ for which $\ox+t_k u_k\in\O$ whenever $k\in\N$. Since weakly bounded sets are norm bounded as well, we get that the sequence $\{u_k\}$ is bounded by norm in $\X$. Thus $\ox+t_k u_k\to\ox$ as $k\to\infty$, which ensures that $\ox+t_k u_k\in\O\cap U$ for all large $k\in\N$, where the neighborhood $U$ is taken from \eqref{Cset}. It tells us by the subamenability of $\O$ that $f(\ox+t_k u_k)\in\Theta$ for large $k$. Employing the assumed differentiability of $f$ at $\ox$ gives us
\begin{equation*}
f(\ox+t_k u_k)= f(\ox)+t_k\nabla f(\ox)u_k+o(t_k)\in\Theta
\end{equation*}
for all $k\in\N$ without loss of generality. Remembering the conic hull definition, we get
\begin{equation*}
\nabla f(\ox)u_k+\frac{o(t_k)}{t_k}\in\frac{\Theta-\oy}{t_k}\subset{\rm cone}(\Theta-\oy)\;\mbox{ with }\;\oy:=f(\ox),\quad k\in\N.
\end{equation*}
Since $\nabla f(\ox)\colon\X\to\Y$ is a bounded linear operator between normed spaces, it is well known to be weak-to-weak continuous; see \cite[Theorem~6.17]{Alip}. Thus $u_k\wto u$ yields $\nabla f(\ox)u_k\wto\nabla f(\ox)u$ as $k\to\infty$. Due to the imposed convexity of $\Theta$ we have the convexity of ${\rm cone}(\Theta-\oy)$, and therefore
\begin{equation*}
\nabla f(\ox)u_k+\frac{o(t_k)}{t_k}\wto\nabla f(\ox)u\in\overline{{\rm cone}\big(\Theta-\oy)}=T_{\Theta}(\oy).
\end{equation*}
Finally, it follows from $\nabla f(\ox)u\in T_{\Theta}(\oy)$ and formula \eqref{tan-msr}, as a consequence of MSQC \eqref{msqcset} by Theorem~\ref{subdchain}, that $u\in T_{\O}(\ox)$. This completes the proof of the claimed assertion in \eqref{wtan}.
\end{proof}\vspace*{-0.05in}

Next we proceed with the study of Dini-Hadamard generalized differential properties of the (always nonsmooth) {\em distance function} that are of their own interest and are needed for our subsequent results. The obtained precise calculations seem to be new in infinite dimensions. We start with the following special case related (not fully) to convex analysis.\vspace*{-0.05in}

\begin{Lemma}[\bf generalized differentiation of the distance function in normed spaces]\label{dist} Let $\O$ be a closed subset of a normed space $\X$, and let $\ph(\cdot):=\dist(\cdot;\O)$ with $\ox\in\O$. If the set $\O$ is convex, then we have the subderivative representation
\begin{equation}\label{disteq}
\d\ph(\ox)(u)=\dist\big(u;T_{\O}(\ox)\big)\;\mbox{ for all }\;u\in\X,
\end{equation}
Furthermore, the fulfillment of \eqref{disteq} for some closed set $\O$ yields the calculation formula
\begin{equation}\label{disteq1}
 \sub^-\ph(\ox)=N^-_{\O}(\ox)\cap\B^*.
\end{equation}
\end{Lemma}\vspace*{-0.1in}
\begin{proof}
We first verify the subderivative representation \eqref{disteq} for closed convex sets $\O$. Pick any $u\in\X$ and readily get for all $t>0$ the equalities
\begin{equation*}
\frac{\ph(\ox+t u)-\ph(\ox)}{t}=\frac{\dist(\ox+tu;\O)}{t}=\dist\Big(u;\frac{\O-\ox}{t}\Big).
\end{equation*}
Since the distance function is Lipschitz continuous, the latter easily implies that
\begin{eqnarray*}
\disp\d\ph(\ox)(u)=\liminf_{t\dn 0}\frac{\ph(\ox+t u)-\ph(\ox)}{t}&=&\liminf_{t \dn 0}\dist\Big(u;\frac{\O-\ox}{t}\Big)\\
&\le&\dist\Big(u;\Limsup_{t\dn 0}\frac{\O-\ox}{t}\Big)=\dist\big(u;T_{\ss\O}(\ox)\big).
\end{eqnarray*}
Employing now the convexity of $\O$, we have for all $t>0$ that $\frac{\O-\ox}{t}\subset T_{\ss\O}(\ox)$, which yields
\begin{equation*}
\dist\big(u;T_{\ss\O}(\ox)\big)\le\liminf_{t\dn 0}\dist\Big(u;\frac{\O-\ox}{t}\Big)=\d\ph(\ox)(u)
\end{equation*}
and thus justifies the subderivative representation \eqref{disteq} in the convex case.\vspace*{0.03in}

To verify further the Dini-Hadamard subdifferential representation \eqref{disteq1}, we are based only on formula \eqref{disteq} without using the convexity of $\O$. Pick any subgradient $v\in\sub^-\ph(\ox)$ and deduce from the duality in \eqref{dhsub} and the subderivative representation \eqref{disteq} that
\begin{equation*}
\la v,u\ra\le\dist\big(u;T_{\O}(\ox)\big)\;\mbox{ whenever }\;u\in\X.
\end{equation*}
This implies, in particular, that $\la v,u\ra\le 0$ for all $u\in T_{\O}(\ox)$, which tells us by \eqref{dual} that $v\in N_{\O}^-(\ox)$. Furthermore, it follows from $0\in T_{\O} (\ox)$ that $\|v\|\le 1$, i.e., $v\in\B^*$. To prove the opposite inclusion in \eqref{disteq1}, take $v\in N_{\O}(\ox)\cap\B^*$ and fix vectors $u\in\X$ and $w\in T_{\O}(\ox)$. Then we have
\begin{equation*}
\la v,u\ra=\la v,u-w\ra+\la v,w\ra\le\|v\|\cdot\|u-w\|\le\|u-w\|.
\end{equation*}
Remembering that $w$ was chosen arbitrarily in $T_{\O}(\ox)$ tells us that $\la v,u\ra\le\dist(u;T_{\O}(\ox))$ by \eqref{disteq}. The latter yields $v\in\sub^-\ph(\ox)$ by \eqref{dhsub}, which verifies \eqref{disteq1} and completes the proof.
\end{proof}\vspace*{-0.05in}

Note that when  the space $\X$ is finite-dimensional, the above subderivative formula \eqref{disteq}---and hence the subdifferential one \eqref{disteq1}---hold for an arbitrary {\em closed} set; see, e.g., \cite[Example~8.53]{rw} with an extensive, nontrivial proof. However, \eqref{disteq} {\em fails} even in Hilbert spaces. Indeed, it is shown in \cite[Theorem~2]{bf} that in {\em any} infinite-dimensional normed space $\X$ there exist a closed set $\O$ and points $\ox,u\in\X$ such that $\d\ph(\ox)(u)<\dist(u;T_\O(\ox))$.\vspace*{0.03in}

Having this in mind, next we prove that both results in Lemma~\ref{dist} hold true when the set $\O\subset\X$ is {\em subamenable} at $\ox$ and the space $\X$ is Banach and reflexive.\vspace*{-0.05in}

\begin{Theorem}[\bf generalized differentiation of the distance function for subamenable set]\label{distam}
Let $\O$ be a nonempty subset of a reflexive Banach space $\X$. Assume that $\O$ is subamenable at $\ox\in\O$ and consider the distance function $\ph(\cdot):=\dist(\cdot;\O)$ associated with $\O$. Then both subderivative and subdifferential formulas \eqref{disteq} and \eqref{disteq1} hold in this case. Furthermore, $\ph$ is epi-differentiable at $\ox$ and its subderivative is represented by
\begin{equation}\label{distam1}
\d\ph(\ox)(u)=\lim_{t\dn 0}\frac{\ph(\ox+t u)-\ph(\ox)}{t}\;\mbox{ for all }\;u\in\X.
\end{equation}
\end{Theorem}\vspace*{-0.2in}
\begin{proof} We first verify the subderivative formula \eqref{disteq}, which readily yields the subdifferential one \eqref{disteq1} by Lemma~\ref{dist}. Observe from the proof of Lemma~\ref{dist} that we always have
\begin{equation*}
\d\ph(\ox)(u =\liminf_{t \dn 0}\dist\Big(u;\frac{\O-\ox}{t}\Big)\le\dist\big(u;T_{\O}(\ox)\big),\quad u\in\X,
\end{equation*}
i.e., the inequality ``$\le$" in \eqref{disteq} is satisfied. To prove the opposite inequality therein in the subamenable case, fix any $\varepsilon>0$ and by \eqref{tan} find sequences $t_k\dn 0$ and $u_k\in\frac{\O-\ox}{t_k}$ such that
\begin{equation*}
\|u-u_k\|-\varepsilon\le\liminf_{t\dn 0}\dist\Big(u;\frac{\O-\ox}{t}\Big)=\d\ph(\ox)(u)<\infty\;\mbox{ for all }\;k\in\N.
\end{equation*}
This tells us that the sequence $\{u_k\}$ is a bounded in the reflexive Banach space $\X$, and so it contains a weakly convergent subsequence. By passing to such a subsequence if needed, suppose without loss of generality that the entire sequence $\{u_k\}$ weakly converges to some vector $\ou\in\X$. Employing now the subamenability of $\O$ at $\ox$ together with $\ox+t_k u_k\in\O$ for all $k\in\N$, we deduce from Proposition~\ref{weaktan} that $\ou\in T_{\O}^w(\ox)=T_{\O}(\ox)$. Since the norm function is weakly lower semicontinuous on $\X$, it gives us
\begin{eqnarray*}
\disp\d\ph(\ox)(u)&=&\liminf_{t\dn 0}\dist\Big(u;\frac{\O-\ox}{t}\Big)\ge\liminf_{k\to\infty}\|u-u_k\|-\varepsilon\\
&\ge&\|u-\ou\|-\varepsilon\ge\dist\big(u;T_{\O}^w(\ox)\big)-\varepsilon=\dist\big(u;T_{\O}(\ox)\big)-\varepsilon,
\end{eqnarray*}
which justify the subderivative formula \eqref{disteq}, and hence the subdifferential one \eqref{disteq1}.\vspace*{0.03in}

It remains to verify the subderivative representation in \eqref{distam1}, which implies the epi-differentiability of $\ph$ at $\ox$ by Lemma~\ref{epilemma}. We get from the above that it is sufficient to show that
\begin{equation}\label{dist-epi}
\d\ph(\ox)(u)=\disp\lim_{t\dn 0}\dist\Big(u;\frac{\O-\ox}{t}\Big)\;\mbox{ for all }\;u\in\X.
\end{equation}
To proceed with proving \eqref{dist-epi}, pick any $u\in\X$ and any $\varepsilon>0$. The obtained formula \eqref{disteq} ensures the existence of a tangent vector $z\in T_{\O}(\ox)$ satisfying
\begin{equation*}
\d\ph(\ox)(u)+\varepsilon=\dist\big(u;T_{\O}(\ox)\big)+\varepsilon\ge\|u- z\|,\quad u\in\X.
\end{equation*}
Since $\O$ subamenable at $\ox$, it is derivable at this point. This gives us a path $z(\cdot)\colon[0,1]\to\X$ such that $z(t)\in\frac{\O-\ox}{t}$ for all $t>0$ with $\lim_{t\dn 0}z(t)=z$. Therefore, for all $u\in\X$ we have
\begin{equation*}
\d\ph(\ox)(u)+\varepsilon\ge\lim_{t\dn 0}\|u-z(t)\|\ge\limsup_{t\dn 0}\dist\Big(u;\frac{\O-\ox}{t}\Big)\ge\liminf_{t\dn 0}\dist\Big(u;\frac{\O-\ox}{t}\Big)=\d\ph(\ox)(u).
\end{equation*}
Since $\varepsilon>0$ was chosen arbitrarily, this verifies \eqref{dist-epi} and finishes the proof.
\end{proof}\vspace*{-0.05in}

The next important result provides a precise formula for calculating the Dini-Hadamard subnormal cone \eqref{dual} to the {\em inverse images}
\begin{equation}\label{inv-im}
\O=f^{-1}(\Th):=\big\{x\in\X\;\big|\;f(x)\in\Th\big\},
\end{equation}
which describe {\em constraint sets} in optimization problems; see Section~\ref{sect06}. Under the {\em metric regularity} and its characterizations for constraint sets, such results are well known for various type of normal cones; see, e.g., the books \cite{bz,c,i,m06,m18,penot,rw}. However, the following theorem derives the normal cone formula under the much weaker metric {\em subregularity} qualification condition \eqref{msqcset}, which ensures the {\em subamenability} of the set $\O$ in \eqref{inv-im}. Moreover, we establish in addition the new {\em bounded multiplier property} for normals to $\O$ in infinite dimensions, a counterpart of which was obtained in \cite[Lemma~2.1]{gm} for the limiting normals from \eqref{lnc} in finite dimensions. The established bounded multiplier property plays a crucial role in deriving new necessary optimality conditions in Sections~\ref{sect06} and \ref{sect07}.

Note also that, although it might be suspected that the subnormal cone formula obtained below is a consequence of the subdifferential chain rule from either Theorem~\ref{lipchainrule} or Theorem~\ref{2chain} with $\th:=\dd_\Th$ therein, this is actually not the case. Indeed, in Theorem~\ref{lipchainrule} the function $\th$ is assumed to be locally Lipschitzian around $\oy$, while the weak$^*$ closedness of the set $\nabla f(\ox)^*N_\Th(\oy)$ required in Theorem~\ref{2chain} may fail even in finite-dimensional spaces.\vspace*{-0.05in}

\begin{Theorem}[\bf normals vectors to inverse images and bounded multiplier property]\label{cahinset}
Let $f\colon\X\to\Y$ be continuously differentiable mapping around $\ox\in\X$, where $\X$ is a reflexive Banach space, while $\Y$ is an arbitrary normed one. Consider the inverse image set $\O$ from \eqref{inv-im}, where $\Theta\subset\Y$ is a closed and convex set with $\oy:=f(\ox)\in\Th$, and assume that MSQC \eqref{msqcset} holds at $\ox$. Then we have the normal cone representation/chain rule
\begin{equation}\label{chainset1}
N_{\O}^-(\ox)=\nabla f(\ox)^*N_{\Theta}(\oy),
\end{equation}
which is strengthened by the following {\sc bounded multiplier property}: for every normal $v\in N_\O^-(\ox)$ there exists a multiplier $\lambda\in N_\Theta(\oy)$ such that
\begin{equation}\label{bounded}
\nabla f(\ox)^*\lambda=v\;\mbox{ and }\;\|\lambda\|\le\kappa\|v\|,
\end{equation}
where the $($uniform$)$ constant $\kappa>0$ is taken from \eqref{msqcset}.
\end{Theorem}\vspace*{-0.15in}
\begin{proof} It is easy to deduce from the definitions that the inclusion $\nabla f(\ox)^*N_\Theta(\oy)\subset N_{\O}(\ox)$ is always satisfied. Thus it remains to verify the bounded multiplier property \eqref{bounded}, which obviously yields the converse inclusion in \eqref{chainset1}, To proceed, consider the Lipschitz continuous functions
\begin{equation*}
\ph(x):=\dist(x;\O)\;\mbox{ and }\;\psi(x):=\kappa\,\dist\big(f(x);\Theta\big),\quad x\in\X.
\end{equation*}
Then the fulfillment of MSQC \eqref{msqcset} at $\ox$ can be rewritten as $\ph(x)\le\psi(x)$ on a neighborhood $U$ of $\ox$. Due to $\ph(\ox)=\psi(\ox)=0$, it is straightforward to deduced from the definitions that
\begin{equation}\label{disteq2}
\d\ph(\ox)(u)\le\d\psi(\ox)(u)\;\mbox{ for all }\;u\in\X,\;\mbox{ and so }\;\sub^-\ph(\ox)\subset\sub^-\psi(\ox).
\end{equation}
Since $\O$ is subamenable at $\ox$, we get by Theorem~\ref{distam} that $\sub^-\ph(\ox)=N_{\O}^-(\ox)\cap\B^*$. Applying further Theorem~\ref{lipchainrule} to the function $\psi$ as the composition of $\kappa\,\dist(\cdot;\Theta)$ and $f$ and then using Lemma~\ref{dist} to calculate the subdifferential of the outer function in the composition yield
\begin{equation*}
\sub^-\psi(\ox)=\kappa\,\nabla f(\ox)^*\big(N_{\Theta}(\oy)\cap\B^*\big).
\end{equation*}
Now take $0\ne v\in N_{\O}^-(\ox)$ and observe from \eqref{disteq2} and the subdifferential formula \eqref{disteq1} of Theorem~\ref{distam} that $\frac{v}{\|v\|}\in\sub^-\ph(\ox)\subset\sub^-\psi(\ox)$. This tells us in turn that there exists $w\in N_{\Theta}(\oy)$ such that $\|w\|\le 1$ and $\frac{v}{\|v\|}=\kappa\,\nabla f(\ox)^*w$. Finally, by setting $\lambda:=\kappa\,\|v\|w$ we conclude that $\lambda\in N_{\Theta}(\oy)$ with $v=\nabla f(\ox)^*\lambda$ and $\|\lambda\|=\kappa\,\|v\|\cdot\|w\|\le\kappa\,\|v\|$, which verifies \eqref{bounded} and thus completes the proof of the theorem.
\end{proof}\vspace*{-0.05in}

It has been well realized in variational analysis that the ``elementary"  generalized differential constructions of the Fr\'echet and Dini-Hadamard types are generally ``nonrobust" meaning that they can be dramatically changed under small perturbations of the initial data. This has a number of negative consequences including poor calculi, etc. To improve the situation, some limiting regularization procedures of type \eqref{lnc} have been developed; see, e.g., \cite{bz,i,m06,m18,penot,rw}. However, implementations of such procedures in infinite dimensions require certain ``normal compactness" assumptions that may be rather restrictive and difficult to verify, especially in non-Lipschitzian settings.\vspace*{0.03in}

The following theorem sheds light on our conceptional understanding the role of {\em subamenability}. It shows that in fairly general infinite-dimensional frameworks (and surely in all finite-dimensional spaces) the Dini-Hadamard and related constructions are indeed {\em robust}, i.e., they are automatically stable with respect to the limiting procedures without imposing any restrictive assumptions.\vspace*{-0.05in}

\begin{Theorem}[\bf robustness of the subnormal cone to subamenable sets]\label{robust} Let $\O$ be a nonempty subset of a reflexive Banach space $\X$. If $\O$ is subamenable at $\ox\in\O$, then we have
\begin{equation}\label{robust1}
N_\O^-(\ox)=N_{\O}^{L-}(\ox):=\big\{v\in\X^*\;\big|\;\exists\,v_k\wsto v,\;\;x_k\st{\O}{\to}\ox\;\;\mbox{with}\;\;v_k\in N_{\O}^{-} (x_k),\;k\in\N\big\}.
\end{equation}
\end{Theorem}\vspace*{-0.15in}
\begin{proof} Since the inclusion ``$\subset$" in \eqref{robust1} is obvious, we need to prove the opposite one. To proceed, pick any $v\in N_\O^{L-}(\ox)$ and find by the definition in \eqref{robust1} such sequences $x_k\st{\O}{\to}\ox$ and $v_k\wsto v$ as $k\to\infty$ that $v_k\in N_\O^-(x_k)$ for all $k\in\N$. By the subamenability of $\O$ at $\ox$, consider its local representation
\eqref{Cset} with some $f\colon\X\to\Y$ and $\Th\subset\Y$ and observe that the imposed MSQC \eqref{msqcset} at $\ox$ holds in fact at all points sufficiently close to $\ox$; in particular, at each $x_k$ with large $k\in\N$. Employing Theorem~\ref{cahinset}, for all such $k$ we find $\lambda_{k}\in N_\Theta(\oy)$ with $\nabla f(x_k)^*\lm_k=v_k$ satisfying the estimate $\|\lambda_k\|\le\kappa\|v_k\|$ {\em uniformly} in $k$. Since the sequence $\{ v_k \}$ is weak$^*$ convergent, the Banach-Steinhaus uniform boundedness principle tells us that the sequence $\{v_k\}$ is norm-bounded in $\X^*$. Thus the latter estimate of $\|\lm_k\|$ ensures that the sequence $\{\lambda_k\}$ is norm-bounded in $\Y^*$. Using further the Banach-Alaoglu theorem in the dual space $\Y^*$ to the normed space $\Y$, we conclude that the sequence $\{\lambda_k\}$ contains a weak$^*$ convergent {\em subnet}. Take such a subnet $\{\lm_{k_\nu}\}$ and find $\lm\in\Y^*$ with $\lm_{k_n\nu}\st{w^*}\to\lm$. Since $f$ from \eqref{Cset} is continuously differentiable around $\ox$, the mapping $\nabla f(\cdot)^*$ is norm-to-norm continuous. This allows us to verify the net weak$^*$ convergence
\begin{equation*}
v_{k_\nu}=\nabla f(x_{k_\nu})^*\lambda_{k_\nu}\wstoo\nabla f(\ox)^*\lambda,
\end{equation*}
which shows therefore that $v=\nabla f(\ox)^*\lambda$.

To justify now that $\lambda\in N_{\Theta}(\oy)$, we rewrite $\lambda_{k_\nu}\in N_\Theta(f(x_{k_\nu}))$ by the convexity of $\Th$ as
\begin{equation*}
\la\lambda_{k_\nu},y-\oy\ra\le 0\;\mbox{ for all }\;y\in\Theta,
\end{equation*}
which implies by passing to the limit along this subnet that $\la\lambda,y-\oy\ra\le 0$ whenever $y\in\Theta$, and hence $\lambda\in N_{\Theta}(\oy)$. Appealing again to Theorem~\ref{cahinset} allows us to conclude that
\begin{equation*}
v=\nabla f(\ox)^*\lambda\in\nabla f(\ox)^*N_\Theta(\oy)=N_\O^-(\ox).
\end{equation*}
This justifies the inclusion ``$\supset$" in \eqref{robust1} and thus completes the proof of the theorem.
\end{proof}\vspace*{-0.05in}

As an immediate consequence of Theorem~\ref{robust} and the discussions in Section~\ref{prelim}, observe that---for subamenable sets $\O$ in finite dimensions---the subnormal cone \eqref{dual} agrees with the limiting normal cone \eqref{lnc} at any point $\ox\in\O$ of the set subamenability.\vspace*{0.03in}

Finally in this section, we relate subamenability with the notion of prox-regularity for functions and sets in Banach spaces, which plays a significant role in variational analysis, particularly in its second-order aspects. Referring the reader to \cite{bert,b,c,m06,pr96,rw} and the bibliographies therein to prox-regularity and its applications in finite and infinite dimensions, let us present now its Dini-Hadamard version in normed spaces. A set $\O\subset\X$ is called {\em prox-regular} at $\ox$ for $\ov\in N_\O^-(\ox)$ if there exist a positive number $r$ together with a neighborhood $U$ of $\ox$ and a bounded neighborhood $V$ of $\ov$ in the norm topology of $\X^*$ such that
\begin{equation}\label{proxregdef}
\la v,u-x\ra\le r\|u-x\|^2\;\mbox{ whenever }\;x,u\in\O\cap U\;\mbox{ and }\;v\in N_\O^-(x)\cap V.
\end{equation}
We also say that a set $\O\subset\X$ is {\em strongly subamenable} at $\ox$ of the mapping $f\colon\X\to\Y$ in its local representation \eqref{Cset} is twice continuously differentiable around $\ox\in\O$ under the validity of the metric subregularity qualification condition \eqref{msqcset}.\vspace*{0.05in}

Now we show that strongly subamenable sets in reflexive Banach space are prox-regular.\vspace*{-0.05in}

\begin{Proposition}[\bf prox-regularity of strongly subamenable sets]\label{proxreg} Let $\O$ be a subset of a reflexive Banach space $\X$. If $\O$ is strongly subamenable at $\ox\in\O$, then it is prox-regular at this point for all the normal vectors $\ov\in N_\O^-(\ox)$.
\end{Proposition}\vspace*{-0.15in}
\begin{proof} Fix any vector $\ov\in N_\O^-(\ox)$ and observe that MSQC \eqref{msqcset} holds around not only at $\ox$, but also---with the same constant $\kappa$---in some neighborhood of $\ox$ that is still denoted by $U$. Let $V$ be a bounded neighborhood of $\ov$ in $X^*$. Then Theorem~\ref{cahinset} tells us that
\begin{equation*}
N_\O^-(x)=\nabla f(x)^*N_\Theta\big(f(x)\big)\;\mbox{ for all }\;x\in\O\cap U,
\end{equation*}
and furthermore we can find a constant $\gg>0$ such that for any $x\in\O\cap U$ and $v\in N_\O^-(x)\cap V$ there exists a multiplier $\lambda\in N_\Theta(f(x))$ with
\begin{equation}\label{sm1}
v=\nabla f(x)^*\lambda\;\mbox{ and }\;\|\lambda\|\le\gg.
\end{equation}
Since $f$ is ${\cal C}^2$-smooth around $\ox$, we  can always shrink $U$ in such a way that $f$ is be represented by
\begin{equation*}
f(u)-f(x)=\nabla f(x)(u-x)+\eta(u,x)\;\mbox{ for all }\;x,u\in U,
\end{equation*}
where $\eta\colon\X^2\to\Y$ satisfies the estimate $\|\eta(u,x)\|\le\beta\|u-x\|^2$ on $U$ with some $\beta>0$. Indeed, the classical mean value theorem allows us to get
\begin{equation*}
\eta(u,x)\colon=\nabla^2f(c_{xu})\big((u-x),(u-x)\big)\;\mbox{ with }\;c_{xu}\in[x,u].
\end{equation*}
Taking now any $x,u\in\O\cap U$ and $v\in N_\O^-(x)\cap V$ with $\lambda$ from \eqref{sm1}, we have
\begin{eqnarray*}
\disp 0&\ge&\langle\lambda,f(y)- f(x)\rangle=\langle\lambda,\nabla f(x)(y-x)\rangle+\langle\lambda,\eta(y,x)\rangle\\
&=&\langle v,y-x\rangle+\langle\lambda,\eta(y,x)\rangle\ge\langle v,y-x\rangle-\gg\beta\|y-x \|^2.
\end{eqnarray*}
Denoting $r:=\gg\beta$, we arrive at \eqref{proxregdef} and thus complete the proof of then proposition.
\end{proof}\vspace*{-0.05in}

As a consequence of the obtained result, we conclude that in the setting of Proposition~\ref{proxreg} the subnormal cone \eqref{dual} agrees with the corresponding {\em proximal normal cone} defined in \cite{bert} in general Banach spaces by following the scheme of \cite{pr96} in finite dimensions.
\vspace*{-0.15in}

\section{Optimality Conditions for General Constrained Problems}\sce\label{sect06}\vspace*{-0.05in}

In this section we develop some applications of the calculus results obtained above to deriving new necessary optimality conditions for general constrained optimization problems in normed spaces. The class of optimization problems under consideration here is described as follows:
\begin{equation}\label{op}
\mbox{minimize }\;\vt(x)\;\mbox{ subject to }\;f(x)\in\Theta
\end{equation}
with $\vt\colon\X\to\oR$, $\Th\subset\Y$, and $f\colon\X\to\Y$ acting between arbitrary normed spaces. We derive below necessary optimality conditions of both primal and dual types by using the Dini-Hadamard generalized differentiation and employing correspondingly the Abadie and metric subregularity qualification conditions, which are now play a role of {\em constraint qualifications}. \vspace*{0.03in}

Let us start with the following {\em primal} necessary optimality condition for local minimizers of \eqref{op} involving the tangent cone \eqref{tan} and subderivative \eqref{subder}. For simplicity we assume the local Lipschitz continuity of the cost function in \eqref{op}, while it follows from the proof that this assumption can be relaxed.\vspace*{-0.05in}

\begin{Theorem}[\bf primal necessary optimality conditions under Abadie constraint qualification]\label{pnoc} Given a local minimizer $\ox$ for \eqref{op} in the normed space setting, assume that $\vt\colon\X\to\oR$ is locally Lipschitzian around $\ox$ and epi-differentiable at this point, that $f\colon\X\to\Y$ is continuously differentiable around $\ox$, and that $\Th\subset\Y$ is locally closed around $\oy:=f(\ox)\in\Th$. Suppose in addition that AQC \eqref{AQC} holds at $\ox$ for $\ph:=\dd_\Th\circ f$. Then we have the implication
\begin{equation}\label{pnoc1}
\big[\nabla f(\ox)u\in T_{\Theta}(\oy)\big]\Longrightarrow\big[\d\vt(\ox)(u)\ge 0\big]\;\mbox{ for all }\;u\in\X.
\end{equation}
\end{Theorem}\vspace*{-0.15in}
\begin{proof} Considering the set of feasible solutions (or the feasible/constraint set) for \eqref{op} given by
\begin{equation}\label{feas}
\O:=\big\{x\in\X\;\big|\;f(x)\in\Th\big\},
\end{equation}
observe that the Abadie qualification condition \eqref{AQC} is written now as
\begin{equation}\label{acq}
T_\O(\ox)=\big\{u\in\X\;\big|\;\nabla f(\ox)u\in T_\Theta(\oy)\big\}.
\end{equation}
It is obvious that the constrained problem \eqref{op} can be equivalently rewritten in the unconstrained form
\begin{equation}\label{op1}
\mbox{minimize }\;\psi(x):=\vt(x)+\delta_\O(x)\;\mbox{ over all }\;x\in\X
\end{equation}
by using the ``infinite penalty." We can directly observe from the subderivative definition \eqref{subder} and the local optimality of $\ox$ in \eqref{op1} that $\d\psi(\ox)(u)\ge 0$ for all $u\in\X$. The summation structure of $\psi$ in \eqref{op1} and the Lipschitz continuity of $\vt$ around $\ox\in{\rm int}(\dom\vt)$ allow us to apply the subderivative sum rule from Theorem~\ref{subdsum}, where the metric inequality \eqref{mscqsum} is obviously satisfied, and hence to arrive at
\begin{eqnarray*}
\d\psi(\ox)(u)=\d\vt(\ox)(u)+\d\delta_\O(\ox)(u)=\d\vt(\ox)(u)+\delta_{\ss T_\O(\ox)}(u)\ge 0\;\mbox{ for all }\;u\in\X.
\end{eqnarray*}
Combining the latter with the Abadie constraint qualification \eqref{acq} tells us that
\begin{equation*}
\nabla f(\ox)u\in T_\Theta(\oy)\Longleftrightarrow u\in T_\O(\ox)\Longrightarrow\d\vt(\ox)(u)\ge 0,
\end{equation*}
which verifies the claimed necessary optimality condition \eqref{pnoc1}.
\end{proof}\vspace*{-0.05in}

It follows from the proof of Theorem~\ref{pnoc} that we have
\begin{equation}\label{feas1}
T_{\O}(\ox)\cap\big\{u\in\X\;\big|\;\d\vt(\ox)(u)<0\big\}=\emp.
\end{equation}
This tells us that the set of {\em improving directions} at $\ox$ does not contain a tangent to the feasible set $\O$ from \eqref{feas} at $\ox$. Both sets in \eqref{feas1} become convex under some additional assumptions, and thus can be separated. In this way, the primal optimality condition of Theorem~\ref{pnoc} yields a certain dual one. The following consequence is a particular realization of this approach, which brings us to an {\em asymptotic} necessary optimality condition expressed of the topological weak$^*$ closure.\vspace*{-0.05in}

\begin{Corollary}[\bf dual asymptotic optimality conditions]\label{dualnoc} Let $\ox\in\O$ is a local minimizer of problem \eqref{op}, where $\O$ is taken from \eqref{feas}. In addition to the assumptions of Theorem~{\rm\ref{pnoc}}, suppose that $\vt$ is Fr\'echet differentiable at $\ox$ and that $\Th$ is convex. Then we have
\begin{equation}\label{dualnoc1}
0\in\nabla\vt(\ox)+\overline{\nabla f(\ox)^*N_\Th(\oy)}^*.
\end{equation}
\end{Corollary}\vspace*{-0.15in}
\begin{proof}
If $\nabla\vt(\ox)=0$, inclusion \eqref{dualnoc1} is trivial. Considering now the case where $\nabla\vt(\ox)\ne 0$ and using condition \eqref{feas1} obtained in Theorem~\ref{pnoc} give us
\begin{equation}\label{feas2}
T_\O(\ox)\cap\big\{u\in\X\;\big|\;\la\nabla\vt(\ox),u\ra<0\big\}=\emp.
\end{equation}
It is easy to see that the tangent cone $T_\O(\ox)$ to the set $\O$ in \eqref{feas} is convex under the assumptions imposed on $f$ and $\Th$. Since the second set in \eqref{feas2} is nonempty and open in $\X$, we employ the convex separation theorem to these sets and find $v\in\X^*\setminus\{0\}$ such that
\begin{equation*}
\la v,z\ra\le\la v,u\ra\;\mbox{ for all }\;z\in T_\O(\ox)\;\mbox{ and }\;u\in\X\;\mbox{ with }\;\la\nabla\vt(\ox),u\ra<0.
\end{equation*}
Since $T_{\O}(\ox)$ is a convex cone and $\nabla\vt(\ox)$ is a continuous linear functional, we deduce from the above that
\begin{equation*}
\la v,z\ra\le 0\le\la v,u\ra\;\mbox{ for all }\;z\in T_\O(\ox)\;\mbox{ and }\;u\in\X\;\mbox{ with }\la\nabla\vt(\ox),u\ra\le 0.
\end{equation*}
The latter clearly implies that $v\in N_\O^-(\ox)$ and $-v\in N_K(0)$, where $K:=\{u\in\X\;|\;\la\nabla\vt(\ox),u\ra\le 0\}$. It is not hard to calculate that $N_K(0)=\{\tau\nabla\vt(\ox)\;|\;\tau\ge 0\}$, which tells us therefore that $-\tau v=\nabla\vt(\ox)$ for some $\tau>0$. Employing finally implication (iv)$\Longrightarrow$(v) in Proposition~\ref{equi}, we conclude that
\begin{equation*}
0=\nabla\vt(\ox)+\tau v\in\nabla\vt(\ox)+N_\O^-(\ox)=\nabla\vt(\ox)+\overline{\nabla f(\ox)^*N_\Th(\oy)}^*,
\end{equation*}
which verifies \eqref{dualnoc1} and completes the proof of the corollary.
\end{proof}\vspace*{-0.05in}

It follows from Example~\ref{ACexample} that the closure operation cannot be generally removed in \eqref{dualnoc1} even in finite dimensions. A sufficient condition for removing the closure therein is that $\Theta$ is a convex polyhedron in $\Y=\R^m$. In particular, if $\Theta:=\{0\}^l\times\R_-^{m-l}$ with some integer $l\in[0,m]$ and if $\X=\R^n$, then \eqref{dualnoc1} reduces to the classical KKT system in nonlinear programming. Observe also that in the case where $\vt$ is continuously differentiable around $\ox$, the feasible set $\O$ from \eqref{feas} is subamenable at $\ox$ provided the fulfillment of the metric subregularity qualification condition \eqref{mscq}. If in addition the space $\X$ is Banach and reflexive, then Theorem~\ref{cahinset} tells us that $N_\O^-(\ox)=\nabla f(\ox)^*N_\Th(\oy)$, which is a weak$^*$ closed subset of $X^*$, and hence the closure operation in \eqref{dualnoc1} can be removed.\vspace*{0.03in}

Now we are going to derive, by using an {\em exact penalization} procedure together with the calculus rules established above, new dual necessary optimality conditions for local minimizers of \eqref{op} without any closure operation as in \eqref{dualnoc1} in the general setting of normed spaces $\X$ and $\Y$ under the {\em metric subregularity constraint qualification}. A strong benefit of metric subregularity is a novel assertion on the {\em boundedness of Lagrange multipliers} with a constructive bound estimate.\vspace*{-0.05in}

\begin{Theorem}[\bf dual necessary optimality conditions with bounded multipliers]\label{dualnopc} Let $\ox$ be a local minimizer of problem \eqref{op} in the general normed space setting. Assume that $\vt$ is Lipschitz continuous around $\ox$ with constant $\ell\ge 0$ and Dini-Hadamard subdifferentially regular at this point, that $f$ is continuously differentiable around $\ox$, that $\Th$ is convex and locally closed around $\oy:=f(\ox)$, and that MSCQ \eqref{msqcset} holds at $\ox$ with some $\kappa>0$. Then there exists $\lambda\in Y^*$ such that
\begin{equation}\label{necessaryop}
0\in\sub^-\vt(\ox)+\nabla f(\ox)^*\lambda,\;\;\lambda\in N_{\Theta}(\oy),\;\mbox{ and }\;\|\lambda\|\le\ell\kappa.
\end{equation}
\end{Theorem}\vspace*{-0.12in}
\begin{proof} First we conclude, by using the (rather elementary) exact penalization procedure from Lemma~4.38 and the arguments on  p.\ 336 of Ioffe's book \cite{i} valid in any normed space, that $\ox$ is a local optimal solution of the unconstrained penalized problem:
\begin{equation}\label{penop}
\mbox{minimize }\;\psi(x):=\vt(x)+\ell\kappa\,\dist\big(f(x);\Theta\big),\quad x\in\X.
\end{equation}
Thus we have by the corresponding subdifferential Fermat rule that $0\in\sub^-\psi(\ox)$. Taking into account the summation structure of the function $\psi$ in \eqref{penop} and the composite structure of the function
\begin{equation*}
\phi(x):=\ell\kappa\,\dist\big(f(x);\Theta\big)=\ell\kappa\,\big(\dist(\cdot;\Theta)\circ f\big)(x),\quad x\in\X,
\end{equation*}
therein, we subsequently apply to $\sub^-\psi(\ox)$ the Dini-Hadamard subdifferential sum rule from Corollary~\ref{subdisum}, the subdifferential representation for the distance function from Lemma~\ref{dist}, and the subdifferential chain rule from Theorem~\ref{lipchainrule}, which all hold under the imposed assumptions in arbitrary normed spaces. In this way we deduce from the aforementioned subdifferential Fermat rule that
\begin{eqnarray*}
\disp
0\in\partial^-\psi(\ox)&=&\sub^-\vt(\ox)+\ell\kappa\,\sub^-\big(\dist(\cdot;\Theta)\circ f\big)(\ox)\\
&=&\sub^-\vt(\ox)+ \ell\kappa\,\nabla f(\ox)^*\big(N_{\Theta}(\oy)\cap\B^*\big)\\
&=&\sub^-\vt(\ox)+\nabla f(\ox)^*\big(N_{\Theta}(\oy)\cap\ell\kappa\B^*\big),
\end{eqnarray*}
which verify \eqref{necessaryop} and completes the proof of the theorem.
\end{proof}\vspace*{-0.12in}

\begin{Remark}[\bf bounded multipliers]\label{bounded multipliers} {\rm The following discussions are now in order:\vspace*{0.03in}

{\bf(i)} Theorem~\ref{dualnopc} not only establishes the KKT-type necessary optimality conditions for local minimizers of constrained optimization problems in normed spaces, by also ensures the dual norm boundedness of Lagrange multipliers for such problems under the metric subregularity constraint qualification. This boundedness, which seems to be not explicitly observed earlier, is strongly required for different issues of variational analysis and optimization; in particular, for the second-order variational theory and applications. Note that the boundedness of multipliers in the necessary optimality conditions \eqref{necessaryop} can be also derived by passing to the infinite penalization \eqref{op1} and employing the calculus results of Theorem~\ref{cahinset}, which is obtained in the framework of reflexive Banach spaces.\vspace*{0.03in}

{\bf(ii)} Since the proof of Theorem~\ref{dualnopc} needs just the corresponding calculus rules in the form of inclusions ``$\subset$" (not the stronger ones as equalities used above), this device works to obtain necessary optimality conditions of type \eqref{necessaryop} for any subdifferential and normal cone notions, which enjoy ``full calculus," i.e., the sum and chain rules with the distance function subdifferentiation, without the assumptions on smoothness and subdifferential regularity. The known constructions of generalized differentiation with such properties are the generalized gradient and normal cone by Clarke \cite{c} in Banach spaces, the sequential limiting subdifferential and normal cone by Mordukhovich \cite{m06} in Asplund spaces, and the approximate $G$-subdifferential and $G$-normal cone by Ioffe \cite{i} in Banach spaces. The latter ones reduce to those in \cite{m06} when the space in question is Asplund and weakly compactly generated (WCG), which includes any separable and any reflexive Banach space.\vspace*{0.03in}

{\bf(iii)} The generalized differential calculi for the constructions discussed in (ii) are based on the completeness of the spaces in question that is needed to employ variational/extremal principles. At the same time, the space completeness is not required in Theorem~\ref{dualnopc}, which allows us to cover, e.g., the spaces $\X={\cal C}^p$, $p\ge 1$, of continuously $p$ times differentiable functions that are particularly important for applications in optimal control. Furthermore, we do not impose any ``normal compactness" assumptions used for limiting procedures in infinite dimensions. Note that the main properties of this type automatic hold when $\dim\Y<\infty$, which is assumed in the recent paper by Guo, Ye and Zhang \cite{gyz}, where necessary optimality conditions (without the boundedness of Lagrange multipliers) are derived for various problems of type \eqref{op} in the case of Banach spaces $\X$ via the generalized differential constructions discussed in (ii). However, the finite dimension of $\Y$  does not allow us to study important classes of optimization problems such as semi-infinite and semidefinite programs, even when $\dim\X<\infty$; see Section~\ref{sect07}.}
\end{Remark}\vspace*{-0.03in}

Finally in this section, observe that in the case where the cost function $\vt$ in \eqref{op} is ${\cal C}^1$-smooth around $\ox$, we can specify the boundedness estimate in \eqref{necessaryop} by putting $\ell=\|\nabla\vt(\ox)\|$ therein.\vspace*{-0.05in}

\begin{Corollary}[\bf enhanced bound estimate for Lagrange multipliers]\label{nablaforlip} In the framework of Theorem~{\rm\ref{dualnopc}}, assume in addition that the cost function $\vt$ is continuously differentiable around $\ox$. Then we have the estimate $\|\lambda\|\le\kappa\|\nabla\vt(\ox)\|$ in \eqref{necessaryop}.
\end{Corollary}\vspace*{-0.15in}
\begin{proof} If $\vt$ is continuously differentiable around $\ox$, then it is clearly Lipschitz continuous around this point, and for any $\ve>0$ the Lipschitz constant of $\vt$ around $\ox$ can be chosen as $\|\nabla\vt(\ox)\|+\ve$. Employing Theorem~\ref{dualnopc}, for each $\ve>0$ we find $\lambda_\ve\in\Y^*$ satisfying \eqref{necessaryop} with $\ell=\|\nabla\vt(\ox)\|+\ve$. It tells us that the net
$\{\lambda_\ve\}_{0<\ve<1}$ is bounded in $\Y^*$. The Banach-Alaoglu theorem allows us to conclude without loss of generality that there exists a multiplier $\lm\in\Y^*$ such that $\lambda_\ve\wstoo\lambda$ as $\ve\dn 0$. It is easy to check by passing to the limit as $\ve\dn 0$ in the KKT conditions
\begin{equation*}
\nabla \vt(\ox)+\nabla f(\ox)^*\lambda_\ve=0\;\mbox{ and }\;\lambda_\ve\in N_\Th(\oy)\;\mbox{ whenever }\;\ve>0
\end{equation*}
that $\nabla\vt(\ox)+\nabla f(\ox)^*\lambda=0$ and $\lambda\in N_\Th(\oy)$. Remembering now that the norm function is weak$^*$ lower semicontinuous in $\Y^*$ and passing to the limit as $\ve\dn 0$ in the estimate
\begin{equation*}
\|\lambda_\ve\|_*\le\kappa\big(\|\nabla\vt(\ox)\|+\ve\big),\quad\ve>0,
\end{equation*}
we arrive at the claimed estimate for $\lm$ and thus complete the proof of the corollary.
\end{proof}\vspace*{-0.25in}

\section{Applications to Semi-Infinite and Semidefinite Programming}\sce\label{sect07}\vspace*{-0.05in}

The concluding section of the paper is devoted to applications of necessary optimality conditions for infinite-dimensional optimization problems of type \eqref{op}, which are obtained in Section~\ref{sect06}, to important classes of {\em semi-infinite programs} (SIPs) and {\em semidefinite programs} (SDPs) formulated via smooth functions in finite-dimensional spaces. Despite the smoothness and finite dimensionality of the initial data, such problems are nonsmooth and infinite-dimensional by their nature, and the results established in Section~\ref{sect06} turn out to be instrumental for the subsequent study of SIPs and SDPs.\vspace*{0.03in}

Starting with inequality-constrained SIPs, we confined ourselves for simplicity to the standard nonlinear model in finite dimensions with compact index sets as studied in the book by Bonnans and Shapiro \cite{bs}. The reader is referred to more recent publications for SIPs with infinite-dimensional decision spaces, nonsmooth functions, and noncompact index sets; see, e.g., \cite{clmp,gl,m18,mn,zy} among other publications with the bibliographies therein. The SIP model under consideration is as follows:
\begin{equation}\label{sip}
\mbox{minimize }\;\vt(x)\;\mbox{ subject to }\;\th(x,s)\le 0\;\mbox{ for all }\;s\in S,
\end{equation}
where $S\ne\emp$ is a {\em compact} metric space of indexes, and where $\vt\colon\R^n\to\R$ and $\th\colon\R^n\times S\to\R$ are continuous. Suppose also, as the {\em standing assumptions} here, that $\vt$ and $\th(\cdot,s)$ are continuously differentiable around the reference points, and that the partial derivative $\nabla_x\th(x,s),$ is continuous in both $(x,s)$.\vspace*{0.03in}

It is easy to see that SIP \eqref{sip} can be reformulated in the infinite-dimensional constrained form \eqref{op} with the same cost function $\vt$ on $\X:=\R^n$, where $\Y:={\cal C}(S)$ is the space of continuous functions on $S$ equipped with the maximum norm, where
\begin{equation}\label{sip-op}
\Th:=\big\{y\in\Y\;\big|\;y(s)\le 0\;\mbox{ for all }\;s\in S\big\},\;\mbox{ and where }\;f(x)(\cdot):=\th(x,\cdot).
\end{equation}
It follows from \cite[Proposition~2.174]{bs} that the mapping $f\colon\R^n\to{\cal C}(S)$ is continuously differentiable around a given feasible point $\ox$ with its Fr\'echet derivative $\nabla f(\ox)\colon\R^n\to{\cal C}(S)$ at $\ox$ calculated by
\begin{equation}\label{part}
\big(\nabla f(\ox)u\big)(s)=\la u,\nabla_{x}\th(\ox,s)\ra\;\mbox{ for all }\;u\in\R^n\;\mbox{ and }\;s\in S.
\end{equation}
Recall also the classical (Riesz) theorem on representing the dual space to ${\cal C}(S)$ as the collection of finite, signed, and regular Borel measures $\mu\colon\mathfrak{B}(S)\to\R$ on the $\sigma$-algebra $\mathfrak{B}(S)$ of Borel subsets of $S$ with
\begin{equation*}
\la\mu,y\ra=\int_{S}y(s)d\mu(s)\;\mbox{ and }\;\|\mu\|=\sup_{E\in\mathfrak{B}(S)}\mu(E)-\inf_{E\in\mathfrak{B}(S)}\mu(E).
\end{equation*}

Among other results of Section~\ref{sect06}, in what follows we consider for brevity applications of only Theorem~\ref{dualnopc} to SIPs and then to SDPs. The next theorem establishes in this way necessary optimality conditions for SIPs \eqref{sip} under consideration, which strictly improves the known results under the smoothness and compactness assumptions made; see Remark~\ref{sipcompare} for more discussions.\vspace*{-0.05in}

\begin{Theorem}[\bf necessary optimality conditions for SIPs with inequality constraints]\label{sipopcon} Let $\ox$ be a local minimizer of SIP \eqref{sip}. In addition to the standing assumptions above, impose that following constraint qualification at $\ox$: there exist a constant $\kappa>0$ and a neighborhood $U$ of $\ox$ such that
\begin{equation}\label{sipcq}
\dist(x;\O)\le\kappa\,\sup_{s\in S}\th^+(x,s)\;\mbox{ for all }\;x\in U,
\end{equation}
where $\O$ is the set of feasible solutions to \eqref{sip}, and where $\th^+(x,s):=\max\{\th(x,s),0\}$. Then we have:
\begin{equation}\label{sipopcon1}
\disp\mbox{there are }\;s_1,\ldots,s_n\in I(\ox)\;\mbox{ and }\;\lambda_1,\ldots,\lambda_n\in\R\;\mbox{ with }\;\left\{\begin{matrix}
\disp\nabla\vt(\ox)+\sum_{i=1}^n\nabla_{x}\lambda_i\th(\ox,s_i)=0,\\\\
\lambda_i\ge 0\;\mbox{ for }\;i=1,\ldots,n,\\\\
\lambda_1+\cdots+\lambda_n\le\kappa\,\|\nabla\vt(\ox)\|,
\end{matrix}\right.
\end{equation}
where $I(\ox):=\{s\in S\;|\;\th(\ox,s)=0\}$ is the set of active indexes associated with $\ox$.
\end{Theorem}\vspace*{-0.15in}
\begin{proof} While reducing SIP \eqref{sip} to the infinite-dimensional optimization problem \eqref{op} with the data in \eqref{sip-op}, let us show that the imposed constraint qualification \eqref{sipcq} is equivalent to the metric subregularity qualification \eqref{msqcset} for the corresponding problem \eqref{op}. To see this, it is sufficient to check that
\begin{equation}\label{ms}
\dist\big(f(x);\Th\big)=\sup_{s\in S}\th^+(x,s)\;\mbox{ for all }\;x\in\R^n,
\end{equation}
where $\Th$ and $f$ are taken from \eqref{sip-op}. Fix $x\in\R^n$ and suppose first that $\sup_{s\in S}\th^+(x,s)=0$. Then $g(x,s)\le 0$ for all $s\in S$, and hence $f(x)\in\Theta$, which verifies \eqref{ms}. If otherwise $\sup_{s\in S}\th^+(x,s)>0$, then the continuity of $\th$ and the compactness of $S$ ensure the existence of $s_x\in S$ with $\sup_{s\in S}\th^+(x,s)=\th(x,s_x)>0$. Taking any
$y(\cdot)\in\Th$, we get that
\begin{equation*}
\|f(x)-y(\cdot)\|=\sup_{s\in S}|\th(x,s)-y(s)|\ge|\th(x,s_x)-y(s_x)|\ge\th(x,s_x)
\end{equation*}
with the equalities holding there when $y(\cdot)=\min\{\th(x,\cdot),0\}$. This tells us that $\dist(f(x);\Th)=\th(x,s_x)$ verifying therefore \eqref{ms} and the fulfillment of MSQC \eqref{msqcset} in \eqref{op} corresponding to SIP \eqref{sip}.\vspace*{0.03in}

Now we are in a position to apply Theorem~\ref{dualnopc} and its supplement in Corollary~\ref{nablaforlip} to problem \eqref{op} corresponding to \eqref{sip}. According to the above discussions, this gives us a Borel measure $\mu\in{\cal C}(S)^*$ with
\begin{equation}\label{sipopcon2}
\begin{array}{ll}
&\nabla\vt(\ox)+\nabla f(\ox)^*\mu =0,\quad\mu\ge 0,\\\\
&\mbox{\rm supp}(\mu)\subset I(\ox),\;\mbox{ and }\;\|\mu\|\le\kappa\|\nabla \vt(\ox)\|.
\end{array}
\end{equation}
It follows from \eqref{part} that we have the representations
\begin{equation*}
\nabla f(\ox)^*\mu=\int_{S}\nabla_{x}\th(\ox,s)d\mu(s):=\bigg(\int_{S}\frac{\sub\th(x,s)}{\sub x_1}d\mu(s),\ldots,\int_{S}\frac{\sub\th(x,s)}{\sub x_n}d\mu(s)\bigg)\in\R^n.
\end{equation*}
Moreover, the calculation in \cite[Example~2.63]{bs} implies that
\begin{equation*}
N_\Th(\oy)=\big\{\mu\in{\cal C}(S)^*\;\big|\;\mu\ge 0\;\mbox{ with }\;{\rm supp}\,\mu \subset\ker f(\ox)\big\}.
\end{equation*}
To proceed further, denote by $\mathfrak{M}$ the collections of measures $\mu\in{\cal C}(S)^*$ satisfying all the conditions in \eqref{sipopcon2}. It is easy to check that the set $\mathfrak{M}$ is convex, weak$^*$ closed, and bounded in ${\cal C}(S)^*$; thus it is weak$^*$ compact there by the Banach-Alaoglu theorem. Then the Krein-Milman theorem tells us that $\mathfrak{M}$ is the weak$^*$ closure of the convex hull of its extreme points. Following the proof of \cite[Lemma~5.110]{bs}, we conclude that any extreme point of $\mathfrak{M}$ has a support with at most $n+1$ points. Since $\mbox{supp}\,\mu\subset I(\ox)$, there exist $s_1,\ldots,s_{n+1}\in I(\ox)$, not necessary different from each other, such that
\begin{equation*}
\mu =\sum_{i=1}^{n+1}\lambda_i\delta(s_i),
\end{equation*}
where $\dd(s)$ stands for the Dirac delta function. Hence it follows from \eqref{sipopcon2} that
\begin{equation}\label{sipopcon3}
\left\{\begin{matrix}
\nabla\vt(\ox)+\disp\sum_{i=1}^{n+1}\nabla_{x}\lambda_i\th(\ox,s_i)=0,\\\\
\lambda_i\ge 0\;\mbox{ for }\;i=1,\ldots,n+1,\\\\
\lambda_1+\cdots+\lambda_{n+1}\le\kappa\,\|\nabla\vt(\ox)\|.
\end{matrix}\right.
\end{equation}
Observe finally that the stationary condition in \eqref{sipopcon3} consists of $n+1$ unknowns $\lambda_i$ and $n$ linear equations. It tells us therefore that we can always find multipliers $\lm_i$ such that only $n$ of them are nonzero. This verifies \eqref{sipopcon1} and thus completes the proof of the theorem.
\end{proof}\vspace*{-0.15in}

\begin{Remark}[\bf comparison with known results]\label{sipcompare} {\rm Let us compare the necessary optimality conditions obtained in Theorem~\ref{sipopcon} with Theorem~5.111 from the book by Bonnans and Shapiro \cite{bs} (see also the references therein), which seems to be---up to now---the strongest known result in this direction for SIPs of type \eqref{sip} under the conventional standing assumptions formulated above. With the exception of the new bounded multiplier assertion $\lambda_1+\ldots+\lambda_n\le\kappa\,\|\nabla\vt(\ox)\|$, the necessary optimality conditions obtained in \cite[Theorem~5.111]{bs} are the same as in \eqref{sipopcon1}, while the main difference is in constraint qualifications: instead of the metric subregularity one \eqref{sipcq}, the previous developments imposed the {\em extended Mangasarian-Fromovitz constraint qualification} formulated as
\begin{equation}\label{rcqforsip}
\la u,\nabla_{x}\th(\ox,s)\ra<0\;\mbox{ for all }\;s\in I(\ox).
\end{equation}
It has been well realized that in the setting under consideration is the same as the fundamental {\em Robinson constraint qualification} \eqref{rob}, which is in turn equivalent to the {\em metric regularity} of the mapping $x\mapsto f(x)-\Th$; see the discussions in Section~\ref{sec03}. On the other hand, the proof of Theorem~\ref{sipopcon} shows that the imposed qualification condition \eqref{sipcq} corresponds to the {\em metric subregularity} of this mapping. The latter condition may be essentially weaker than \eqref{rcqforsip} as demonstrated by simple examples. Consider, e.g., the following one-dimensional SIP \eqref{sip} with linear constraints (the cost function can be arbitrary, since it doesn't enter the constraint qualifications): $\th(x,s)=sx$ with $S:=[0,1]$. It is obvious that \eqref{rcqforsip} fails at $\ox=0$, while \eqref{sipcq} holds due to the polyhedrality of the feasible set $\O:=\{x\in\R\;|\;x\le 0\}$.}
\end{Remark}\vspace*{-0.03in}

Another essential advantage of the metric subregularity constraint qualification in comparison with \eqref{rcqforsip} is that we can easily use \eqref{sipcq} to handle {\em equality constraints} by reducing them to two inequalities. This trick obviously does not work for \eqref{rcqforsip}, since it brings us to the triviality in necessary conditions.\vspace*{0.05in}

Now we derive in this way new necessary optimality conditions for the following class of SIPs with both inequality and (for simplicity just one) equality constraints:
\begin{eqnarray}\label{sipeq}
\disp
&\mbox{minimize}&\vt(x)\;\mbox{ for }\;x\in\R^n\\\nonumber
&\mbox{subject to}&\th(x,s)\le 0\;\mbox{ for }\;s\in S,\\\nonumber
&&\psi(x,t)=0\;\mbox{ for }\;t\in T,
\end{eqnarray}
where $S$ and $T$ are compact metric spaces, and where the functions $\vt$ and $\th$ are taken from \eqref{sip}, while $\psi(x,t)$ satisfies the same standing assumptions as $\th(x,s)$. Here is a consequence of Theorem~\ref{sipopcon}.\vspace*{-0.05in}

\begin{Corollary}[\bf SIPs with equality and inequality constraints]\label{sipeqop} Let $\ox$ be a local minimizer of SIP \eqref{sipeq}. In addition to the standing assumptions, impose the following metric subregularity constraint qualification at $\ox$: there exist a number $\kappa>0$ and a neighborhood $U$ of $\ox$ such that
\begin{equation}\label{sipeqcq}
\dist(x;\O)\le\kappa\,\sup_{s\in S}\th^+(x,s)+\kappa\,\sup_{t\in T}|\psi(x,t)|\;\mbox{ whenever }\;x\in U,
\end{equation}
where $\O$ is the set of feasible solutions to \eqref{sipeq}. Then for some pairs of indexes $(s_1,t_1),\ldots,(s_n,t_n)\in S\times T$ and Lagrange multipliers $(\lambda_1,\mu_1),\ldots,(\lambda_n,\mu_n)\in\R^2$ we have the optimality conditions
\begin{equation}\label{sipeqop1}
\left\{\begin{matrix}
\disp\nabla\vt(\ox)+\sum_{i=1}^n\lambda_i\nabla_x\th(\ox,s_i)+\sum_{i=1}^n\mu_i\nabla_x\psi(\ox,t_i)=0,\\\\
\lambda_i\ge 0\;\mbox{ and }\;\lambda_i\th(\ox,s_i)=0\;\mbox{ for }\;i=1,\ldots,n,\\\\
\disp\sum_{i=1}^n\lambda_i+\sum_{i=1}^n|\mu_{i}|\le 2\kappa\,\|\nabla\vt(\ox)\|.
\end{matrix}\right.
\end{equation}
\end{Corollary}\vspace*{-0.15in}
\begin{proof} The equality constraint in \eqref{sipeq} is obviously represented by the two inequalities
\begin{equation*}
\psi(x,t)\le 0\;\mbox{ and }\;-\psi(x,t)\le 0\;\mbox{ for all }\;t\in T.
\end{equation*}
Thus now we are in the framework of SIP \eqref{sip} with only the inequality constraints. To use Theorem~\ref{sipopcon} in this framework, we need to show that the assumed constraint qualification \eqref{sipeqcq} implies the fulfillment of \eqref{sipcq} for the extended inequality system. Fix $x\in\R^n$ and observe to this end that
\begin{equation*}
\sup_{s\in S}\th^+(x,s)+\sup_{t\in T}|\psi(x,t)|\le 2\sup_{(s,t)\in S\times T}\big(\max\big\{\th^+(x,s),(-\psi)^+(x,t),\psi^+(x,t)\big\}\big).
\end{equation*}
This tells us that the imposed qualification condition \eqref{sipeqcq} yields
\begin{equation*}
\dist(x;\O)\le2\kappa\sup_{(s,t)\in S\times T}\big(\max\big\{\th^+(x,s),(-\psi)^+(x,t),\psi^+(x,t)\big\}\big),
\end{equation*}
which is the one we need in \eqref{sipcq}. Applying Theorem~\ref{sipopcon} to the latter SIP verifies \eqref{sipeqcq}.
\end{proof}

Finally in this section, we consider problems of {\em semidefinite programming} (SDP) the importance of which has been highly recognized not only in continuous optimization but also in a variety of applications coming from optimal control, combinatorial optimization, statistics, etc.; see, e.g., \cite{bs} and the references therein. To formulate this class of optimization problems, denote by ${\cal S}^m$ the space of $m\times m$ symmetric matrix endowed with the usual Frobenius norm and recall that the notation $A\preceq 0$ indicates that the symmetric matrix $A$ is negative semidefinite. Then the basic SDP is given by:
\begin{eqnarray}\label{sdp}
\mbox{minimize}\;\vt(x)\;\mbox{ subject to }\;\Phi(x)\preceq 0\;\mbox{ and }\;\Psi(x)=0,
\end{eqnarray}
where $\vt\colon\R^n\to\R$, $\Phi\colon\R^n\to{\cal S}^m$, and $\Psi\colon\R^n\to{\cal S}^m$ are continuously differentiable around the point in question. Given $A\in{\cal S}^m$, denote by $\sigma(A)$ and by  $\|A\|_{\max}$ the largest eigenvalue of $A$ and the largest absolute value entries of $A$, respectively. It is well known that
\begin{equation*}
\sigma(A)=\sup\big\{\la s,As\ra\big|\;\|s\|=1\big\}\;\;\mbox{ for any }\;A\in{\cal S}^m.
\end{equation*}

The following theorem provides necessary optimality condition for SDP \eqref{sdp} by reducing it to the corresponding SIP and employing the results established above.\vspace*{-0.03in}

\begin{Theorem}[\bf necessary optimality conditions for SDPs]\label{sdpkkt} Let $\ox $ be a local minimizer of SDP \eqref{sdp}. In addition to the standing assumptions, impose the constraint qualification at $\ox$ of the metric subregularity type: there exist a number $\kappa>0$ and a neighborhood $U$ of $\ox$ such that
\begin{equation*}
\dist(x;\O)\le\kappa\,\big(\sigma^+(\Phi(x))+\|\Psi(x)\|_{max}\big)\;\mbox{ for all }\;x\in U,
\end{equation*}
where $\O$ is the set of feasible solutions to \eqref{sdp}, and where $\sigma^+(A):= \max\{0,\sigma(A)\}$. Then we can find vectors $s_1,\ldots,s_n\in\R^m$ and real numbers $\lambda_1,\ldots,\lambda_n$, $\mu_{11},\ldots,\mu_{mm}$ satisfying the conditions
\begin{equation}\label{sdpop1}
\left\{\begin{matrix}
\disp\nabla\vt(\ox)+\sum_{i=1}^n\lambda_i\nabla\Phi(\ox)(s_i,s_i)+\sum_{i=1}^m\sum_{j=1}^m\mu_{ij}\nabla\Psi_{ij}(\ox)=0,\\\\
\|s_i\|=1,\quad\lambda_i\ge 0,\quad\lambda_i\la s_i,\Phi(\ox)s_i\ra=0\;\mbox{ for }\;i=1,\ldots,n,\\\\
\disp\sum_{i=1}^n\lambda_i+\sum_{i=1}^m\sum_{j=1}^m|\mu_{ij}|\le 2\kappa\,\|\nabla\vt(\ox)\|,
\end{matrix}\right.
\end{equation}
where $\nabla\Phi(\ox)(s_i,s_i)$ is a vector in $\R^n$ having the $j^{th}$ component $\big\la s_i,\frac{\partial_j\Phi}{\partial_{x_j}}(\ox)s_i\big\ra$.
\end{Theorem}\vspace*{-0.15in}
\begin{proof} To reduce SDP \eqref{sdp} to the SIP form \eqref{sipeq}, define
\begin{equation}\label{sdp1}
\begin{array}{ll}
\th(x,s):=\la s,\Phi(x)s\ra\;\mbox{ for }\;s\in S\;\mbox{ with }\;S:=\big\{s\in\R^m\big|\;\|s\|=1\big\},\\\\
\psi(x,t):=\Psi_t(x)\;\mbox{ for }\;t\in T\;\mbox{ with }\;T:=\big\{(i,j)\big|\;i,j=1,\ldots,m\big\},
\end{array}
\end{equation}
where $\Psi_t(x)$ denotes the $t=(i,j)$-entry of the matrix $\Psi(x)$. Then it is easy to see that the given SDP \eqref{sdp} reduces the SIP formulation \eqref{sipeq}, where the equality constraint function $\psi(x,t)$ is generally vector-valued. Furthermore, observe that for all $x\in\R^n$ we have
\begin{equation*}
\sup_{s\in S}\th(x,s)=\sigma^+\big(\Phi(x)\big)\;\mbox{ and }\;\max_{t\in T}|\psi(x,t)|=\|\Psi(x)\|_{\max}.
\end{equation*}
The latter tells us that the SIP constraint qualification condition \eqref{sipcq} holds at $\ox$. Using now Corollary~\ref{sipeqop}, we arrive at the KKT type conditions in \eqref{sdpop1} with at most $n$ nonzero scalars $\mu_{ij}$.
\end{proof}\vspace*{-0.15in}

\section{Concluding Remarks}\label{conc}\sce\vspace*{-0.05in}

The main thrust of this paper is to show that appropriate tools of variational analysis and generalized differentiation allow us to work in incomplete normed spaces without using limiting procedures and rather restrictive assumptions. We exploit for these purposes primal and dual generalized differential constructions for extended-real-valued functions on normed spaces that are revolved around subderivatives and Dini-Hadamard subgradients. Major calculus rules are developed for both subderivatives and subgradients under the weakest qualification conditions; namely, the metric subregularity and Abadie ones. Based on the established calculus rules, we derive necessary optimality conditions for a broad class of constrained problems in nonsmooth optimization. The obtained results are largely new even in finite dimensions. Furthermore, the developed variational machinery leads us to novel necessary optimality conditions for problems of semi-infinite and semidefinite programming with smooth data in finite-dimensional spaces by reducing them to the nonsmooth infinite-dimensional framework of this paper.

In our future research, we plan to develop news applications of the approach and results obtained here to remarkable classes of problems from the calculous of variations and optimal control that are naturally formulated in incomplete normed spaces.
\end{document}